\documentclass[preprint, 12pt, numbered]{elsarticle}

\usepackage{color}
\usepackage{a4wide}
\usepackage{amssymb}
\usepackage{amsmath}
\usepackage{graphicx}
\usepackage{hyperref}
\usepackage{mathrsfs}
\usepackage{caption}
\usepackage{subcaption}
\usepackage{textcomp}
\usepackage{multirow}
\usepackage{array}
\usepackage{enumitem}
\hypersetup{
  colorlinks,
  citecolor=blue,
  filecolor=green,
  linkcolor=cyan,
  urlcolor=purple
}
\begin{document}
\begin{frontmatter}
\title{Addressing the issue of mass conservation error and the connected problem of Carbuncle formation}

\author[add1]{Vinnakota Mythreya\corref{cor1}\fnref{fn1}}
\cortext[cor1]{Corresponding Author}
\ead{myth.vinna007@gmail.com}
\fntext[fn1]{Doctoral Candidate}

\author[add1]{S R Siva Prasad Kochi\fnref{fn1}}
\ead{siva.ksr@gmail.com}

\author[add1]{ M. Ramakrishna\fnref{fn2}}
\ead{krishna@ae.iitm.ac.in}
\fntext[fn2]{Professor}

\address[add1]{Department of Aerospace Engineering, IIT Madras}

\begin{abstract}
  We study mass conservation errors (momentum density spike) and the related phenomenon of post shock oscillations in numerical solutions of compressible Euler equations. These phenomena
  and their causes have been reported in literature \cite{roberts1990,arora1997}. In this paper, first, we compare the mass conservation and
  post shock oscillation errors obtained using combinations of different numerical methods (Finite Volume, Finite Difference with WENO and DG with simple WENO limiter)
  and upwind flux functions (ROE, AUSM\textsuperscript{+}-up, and others) for moving shocks, modelled using one-dimensional Euler equations. Next, the mass conservation error
  is quantified for stationary shocks modelled using one-dimensional, quasi-one-dimensional and two-dimensional Euler equations. It is shown that
  using a fine mesh or refining mesh near shocks using multiple over set meshes lead to mitigation of the mass conservation error. This is demonstrated using the 
  problem of flow through a variable area duct modelled using quasi-one-dimensional Euler equations.  The link between mass conservation error and carbuncle formation is shown and preliminary results indicating that the carbuncle can be cured using multiple overset
  meshes are also shown.

  \end{abstract}
\end{frontmatter}
\section{Introduction} \label{Sec:Introduction}
For computing numerical solutions of the compressible Euler equations, upwind biasing of spatial derivatives is done. To achieve this upwind biasing, various
techniques are used, two of which are using approximate Riemann solvers \cite{roe1981} and Flux vector splitting \cite{steger1981}. These techniques, in addition to 
accounting for the direction of propagation of waves, often (not always) introduce ``natural dissipation'' and make the numerical methods stable. However, in the
presence of stationary or moving shocks, this dissipation leads to smearing of shocks, post shock oscillations, and other errors in the numerical solution.

The post shock oscillations error and mass conservation error or momentum spike error, have been reported in a number of papers in literature. 
For slowly moving shocks, Roberts \cite{roberts1990}, investigated the source of error and evaluated the performance of Roe and Osher flux with respect to these errors.
Jin et al \cite{jin1996} explained the source of the post shock oscillation errors from the perspective of smearing due to numerical viscosity. Lin \cite{lin1995} proposed
a modification to
the ROE flux to suppress the post shock oscillation and also reported that the error in the solution is also dependent of direction of motion of shock.
Arora et al \cite{arora1997}, explained the cause of the post shock oscillations, remarked that this problem may be unavoidable for shock capturing schemes without
significant increase in computational effort and suggested some ways to overcome these problems.
Xu \cite{xu1999} asked the question `Does Perfect Riemann Solver Exist'. Xu analysed the dissipation mechanism
in the Godunov scheme, consisting of the `gas evolution stage' `for numerical fluxes across a cell interface' and the `projection stage' `for the reconstruction of constant
state inside each cell'. Xu remarked that the numerical dissipation is solely provided by the `projection stage' and to compute numerical solutions with discontinuities,
addition of explicit dissipation is needed in the `gas evolution stage'. We refer to \cite{stiriba2003, johnsen2013, kitamura2019} for recent analysis of the post shock 
oscillations problem and other `shock anomalies'.

Momentum spike or mass flux errors in steady state numerical solutions with shocks is another such shock anomaly. Barth \cite{barth1989} reported that momentum spike with
error as high as 40\% can occur in numerical solutions. 
When the flux splitting or approximate Riemann solver used in the numerical method leads to smearing of shock, generally an error in mass conservation equation
is introduced. Jin et al \cite{jin1996} show that this is similar to adding dissipation in the mass conservation equation.

In this paper we study the post shock oscillation error and mass conservation error. We define invariants across a moving shock, modelled using one-dimensional Euler
equations and quantify the error in the numerical solution based on the invariant associated with mass conservation.
We compare the performance of different numerical flux functions, namely, the ROE Flux \cite{roe1981}, ROE Flux with Harten Hyman 2 Fix \cite[p. 266]{harten1983},
Osher's Flux with P-Ordering \cite[Section 12.3.1, p. 393]{toro2013}, Osher's Flux with O-Ordering \cite[Section 12.3.2, p. 397]{toro2013}, AUSM\textsuperscript{+}-up
flux \cite{liou2006} and the global Lax-Freidrichs flux.
We compare the errors in
the numerical solutions obtained using these numerical flux functions and numerical methods with different formal order of accuracies. We show that there is no one particular
flux function that consistently performs better than others for different shock speeds, shock movement directions and order of accuracies of the numerical methods.

Later, for steady state solutions with shocks, we demonstrate
the error in mass conservation or `mass conservation error' or `mass flux error', with the help of test problems having shocks in solution,
for one-dimensional Euler equations, quasi-one-dimensional Euler equations and two-dimensional Euler equations. We show how this mass flux error
varies with respect to different parameters like the formal order of accuracy of the numerical method, etc.

It is known that for the one-dimensional Euler equations, using less dissipative fluxes like the Roe flux can lead to capturing a normal shock without smearing and mass
conservation error.
This changes with the introduction of viscous fluxes. Even at large Reynolds numbers and comparatively low magnitude of viscous fluxes, the shock gets dissipated and this
coupled with the ROE flux leads to significant errors in mass conservation. We demonstrate this for the one-dimensional viscous fluid flow equations (Newtonian, Navier-Stokes)
for the problem of normal shock. We show that one way to mitigate this error is using a sufficiently fine mesh to resolve the shock.
We indicate an efficient way of resolving flow near a shock is using multiple overset meshes, using which shocks can be captured with
less error and demonstrate this using the problem of flow through a variable area duct. We also show the connection between mass conservation error and carbuncle formation, using the problem of flow over a cylinder.
We show a way to cure the carbuncle by refining the mesh near the shock using multiple overset mesh and present preliminary results.

In this paper, we use three numerical methods, namely, the finite volume method, the Shu-Osher conservative finite difference scheme and Discontinuous Galerkin method.
For high order finite
volume and finite difference methods, we use high-order WENO \cite{jiang1996} or linear reconstruction.
The mass conservation error and post shock oscillations seem essentially independent of the underlying numerical methods. A comparison of the slight differences among the solutions obtained using these methods is presented, wherever they are interesting.
For discretisation of time derivatives, we use the TVD-RK3 method.

The rest of the paper will be organised as follows: In section \ref{Sec:NumMethod}, we give a brief description of the finite volume method, the Shu-Osher
conservative finite difference method, high-order WENO and linear reconstructions, the Discontinuous Galerkin method with simple WENO limiter, TVD-RK3 method,
different numerical fluxes and flux splittings. We give labels for the different numerical methods used. In section \ref{Sec:Verification}, the implementation
of WENO and DG schemes are verified
using test problems of Burgers equation with source term and the Isentropic Euler vortex. In section \ref{Sec:postShockOscillations} we compare the performance of different
numerical flux functions and numerical methods of different order of accuracy for moving shock problems. We show that there are certain problem parameters for which
the Roe flux produces lesser error than the Osher flux. We underscore the importance of characteristic-wise reconstruction by giving examples of problems for which
doing component-wise reconstruction leads to `NAN's in the computations.
In section \ref{Sec:massFluxErr} we show how the mass flux error varies in numerical solutions having stationary shocks,
for different problems and different numerical methods. We explain the cause of the mass flux error, indicate technique to mitigate it and demonstrate it by
applying the technique for two problems. In section \ref{Sec:refinementAndOverset}, we show that refining near the shock using multiple overset meshes near shocks leads to the mitigation of mass flux error and demonstrate this by using two levels of
over set meshes ( One mesh, overset with a finer mesh, which is in turn overset with a finer mesh ), for the problem of flow through a variable area duct with a normal
shock. We also show the link between the mass conservation error and carbuncle formation and that it can be cured by refining near the shock using overset mesh with two levels of refinement.
We end the paper with concluding remarks in section \ref{Sec:Conclusions}.


\section{ Numerical methods }\label{Sec:NumMethod}
%
%
%
%
%
In this section the numerical methods used, namely, the finite volume and the Shu-Osher conservative finite difference method with Weighted Essentially Non-oscillatory
(WENO) reconstruction and Discontinuous Galerkin method are described. Methods used for time discretisation, namely Total Variation Diminishing-Three stage Runge
Kutta (TVD-RK3) and Butcher's six stage Runge Kutta time discretisation are described. Also, the flux splitting and approximate Riemann solvers used will be
described shortly. We start with the description of the finite volume scheme.

\subsection{Finite Volume method}
Consider a hyperbolic conservation law of the form
\begin{equation}
  \label{eq:hypConsLaw}
  \frac{\partial }{\partial t}Q(x,t) + \frac{\partial}{\partial x}E(Q(x, t)) = 0
\end{equation}
The physical domain is divided into `n' cells, with the $i^{th}$ cell having size equal to $\Delta x$. Integrating equation \ref{eq:hypConsLaw}, over the $i^{th}$ cell with
boundaries $[x_{i-\frac{1}{2}}, x_{i+\frac{1}{2}}]$ between times $t_n$ and $t_{n+1}$ (separated by $\Delta t$ ), we get
\begin{equation}
  \label{eq:FVEqn}
  (\bar{Q}_i^{n+1} - \bar{Q}_i^{n})\Delta x + \int\limits_{t_n}^{t_{n+1}}\left(E(x_{i+\frac{1}{2}}, t) - E(x_{i-\frac{1}{2}}, t)\right)dt = 0,
\end{equation}
where $\bar{Q}^n_i$ is the cell average of Q over the $i^{th}$ cell at time level $t_n$. The time integral of the flux ($E$) is approximated using a numerical flux
function ($\hat{E}$)
\begin{equation}
  \label{eq:FluxFunction}
  \int\limits_{t_n}^{t_{n+1}}E(x_{i+\frac{1}{2}}, t)dt = \hat{E}(Ql^{n}_{i+\frac{1}{2}}, Qr^{n}_{i+\frac{1}{2}})\Delta t,
\end{equation}
where $Ql^{n}_{i+\frac{1}{2}}, Qr^{n}_{i+\frac{1}{2}}$ are left and right biased approximations to $Q(x_{i+\frac{1}{2}}, t_n)$. These approximations will be obtained using
different high order reconstruction procedures that will be described later. The flux function can be based on the approximate Riemann solver of ROE or the AUSM splitting and
others that will be described later. Next we briefly describe the Shu-Osher conservative finite difference method.

\subsection{Shu-Osher Conservative finite difference scheme}\label{Sec:shuOsherConsFinDiffSch}
Let the computational domain consist of grid points uniformly spaced in the physical domain, with grid point spacing equal to $\Delta x$. A function $h(x, t)$ is defined
such that the sliding average of $h(x, t)$ over a length $\Delta x$ is equal to $E(x, t)$, that is,
\begin{equation}
  \label{eq:ImpFn}
  \frac{1}{\Delta x}\int\limits_{-\frac{\Delta x}{2}}^{\frac{\Delta x}{2}}h(x+y, t)dy = E(x, t)
\end{equation}
Taking a partial derivative of equation~(\ref{eq:ImpFn}) with $x$, we get
\begin{equation}
  \label{eq:derivInTermsOfImpFn}
  \frac{\partial E}{\partial x}\bigg|_{x=x_o} = \frac{h(x_o+\frac{\Delta x}{2}, t) - h(x_o-\frac{\Delta x}{2}, t)}{\Delta x}
\end{equation}
We refer to Barry Merriman \cite{merriman2003} for detailed explanation and analysis of the Shu-Osher conservative finite difference scheme.

Using the method of lines and equations~(\ref{eq:hypConsLaw}), and (\ref{eq:derivInTermsOfImpFn}), a semi-discrete form of equation~(\ref{eq:hypConsLaw}) is obtained at
$x=x_o, t=t_o$,
which is
\begin{equation}
  \label{eq:discreteHypConsLaw}
  \frac{\partial Q}{\partial t}\bigg|_{x=x_0, t=t_o} + \frac{h(x_o+\frac{\Delta x}{2}, t_o) - h(x_o-\frac{\Delta x}{2}, t_o)}{\Delta x} = 0
\end{equation}

\subsection{Upwinding and Flux-Splitting}\label{Sec:Upwinding}
To account for propagation along the characteristic directions,
upwind biasing of spatial derivatives is needed. This can be achieved by using flux splitting and appropriate biasing of approximations involving the split fluxes.
\begin{equation}
  \label{eq:FSplitting}
  E^{\pm} = \frac{1}{2}\left(E(Q) \pm \hat{A} Q\right).
\end{equation}
Equation (\ref{eq:FSplitting}) gives the expression for the split fluxes.
Different choices of $\hat{A}$ leads to different flux splittings which will be described in detail later.
%

The semi-discrete form of the hyperbolic conservation law incorporating flux
splitting becomes 
\begin{equation}
  \label{eq:discreteHypConsLawFS}
  \frac{\partial Q}{\partial t}\bigg|_{x=x_0, t=t_o} + \frac{h^{+}(x_o+\frac{\Delta x}{2}, t_o) - h^{+}(x_o-\frac{\Delta x}{2}, t_o)}{\Delta x} + \frac{h^{-}(x_o+\frac{\Delta x}{2}, t_o) - h^{-}(x_o-\frac{\Delta x}{2}, t_o)}{\Delta x}= 0,
\end{equation}
where
\begin{equation}
  \label{eq:ImpFnFS}
  \frac{1}{\Delta x}\int\limits_{-\frac{\Delta x}{2}}^{\frac{\Delta x}{2}}h^{\pm}(x+y, t)dy = E^{\pm}(x, t).
\end{equation}
A high-order linear reconstruction or WENO reconstruction procedure is used to obtain approximations to $h^{+}$ and $h^{-}$ using left and right biased stencils,
respectively. It is described next.
\subsection{Linear reconstruction and WENO reconstruction procedures}\label{Sec:WENOPROC}
Using high-order reconstruction (linear reconstruction) to approximate $h^{+}$ and $h^{-}$ in the presence of discontinuities will lead to oscillations in the
solution. To avoid this, we use WENO reconstruction wherever necessary.

WENO reconstruction was introduced by Liu, Osher and Chan in 1994 \cite{liu1994}.
Jiang et al gave a framework to build high order(of order $2r - 1$ for $r=2, 3, ...$) WENO schemes \cite{jiang1996}. Changes to these schemes were 
proposed \cite{borges2008, castro2011, wu2016} to avoid loss of accuracy near critical points. Two such schemes are, WENO-Z (or ZWENO) proposed by Borges et al
\cite{borges2008, castro2011} and WENO-NP3 \cite{wu2016} proposed by Wu et al.

Equation~(\ref{eq:discreteHypConsLawFS}) is used to advance from time $t_n$ to $t_{n+1}$.  At grid point
$x_i$, approximations $\hat{h}^{\pm}_{i+\frac{1}{2}}$ and $\hat{h}^{\pm}_{i-\frac{1}{2}}$ (subscript $n$, indicating time level, is dropped for brevity) to
$h^{\pm}(x_{i+\frac{1}{2}}, t_n)$ and $h^{\pm}(x_{i-\frac{1}{2}}, t_n)$, respectively, are needed. These, for a $2r-1$ reconstruction are given by the following equations:
\begin{align}
  \label{eq:fluxReconsApprox}
  \hat{h}^{\pm}_{i+\frac{1}{2}} &= \sum\limits_{j=1}^r \omega^{\pm}_j H^{\pm}_j,~
  \omega^{\pm}_j = \frac{\tilde{\omega}_j}{\bar{\omega}^{\pm}},~
  \bar{\omega}^{\pm} = \sum\limits_{j=1}^r \tilde{\omega}^{\pm}_j\\
  \label{eq:omegaWenoZ}
  \tilde{\omega}^{\pm}_j &= \gamma^{\pm}_j\bigg(1+\bigg(\frac{\tau^{\pm}}{\beta^{\pm}_j + \epsilon}\bigg)^p\bigg),\text{ for ZWENO (ZW) reconstruction \cite{borges2008}, and}\\
  \label{eq:linHiOrd}
  \tilde{\omega}^{\pm}_j &= \gamma^{\pm}_j\text{ for linear reconstruction (LR)}.
\end{align}
where $\tau^{\pm}$ are high order smoothness indicators \cite{borges2008}, $\beta^{\pm}_j$ are the Jiang-Shu smoothness indicators \cite{jiang1996}, $\gamma^{\pm}_j$ are linear
weights, $\epsilon = 10^{-14}$, $p = r-1$ (unless specified otherwise).
$\gamma^{\pm}_j \text{ and } H^{\pm}_j$ for a third order reconstruction
are given by:
\begin{align}
  \label{eq:gammasThirdOrd}
  \gamma^{+}_1 = \frac{1}{3},~ \gamma^{+}_2 = \frac{2}{3}, &~ 
  \gamma^{-}_1 = \frac{2}{3},~ \gamma^{-}_2 = \frac{1}{3}.\\
  H^{+}_1 = \frac{3E^{+}_i - E^{+}_{i-1}}{2},~  H^{+}_2 = \frac{E^{+}_i + E^{+}_{i+1}}{2},&~
  H^{-}_1 = \frac{E^{-}_i + E^{-}_{i+1}}{2},~  H^{-}_2 = \frac{3E^{-}_{i+1} - E^{-}_{i+2}}{2}.
\end{align}
Formulae for $\beta^{\pm}_j$ for WENO-NP3 or ZWENO3 ($r=2$) reconstruction, for which a stencil of 3 points is used
are given below.
\begingroup
\allowdisplaybreaks
\begin{align}
  \beta^{+}_1 = (E^+_{i-1} - E^+_i)^2, \beta^{+}_2 = (E^+_{i+1} - E^+_i)^2,&
  \beta^{-}_1 = (E^-_{i+1} - E^-_i)^2, \beta^{-}_2 = (E^-_{i+1} - E^-_{i+2})^2,\\
  \dot{\beta}^{+} = \frac{1}{4}(E^+_{i-1} - E^+_{i+1})^2 + &\frac{13}{12}(E^+_{i-1} - 2E^+_{i} + E^+_{i+1})^2,\\
  \dot{\beta}^{-} = \frac{1}{4}(E^-_{i} - E^-_{i+2})^2 + &\frac{13}{12}(E^-_{i} - 2E^-_{i+1} + E^-_{i+2})^2,\\
  \tau^{\pm} = \tau^{\pm}_{NP} = \Bigg| \dot{\beta}^{\pm} - \frac{\beta^{\pm}_1 + \beta^{\pm}_2}{2}\Bigg|^{1.5},~
  E^{\pm}_k &= E^{\pm}(x_k, t_n)\text{, for } k=i-1, i, i+1, i+2.
\end{align}
\endgroup
We refer to \cite{borges2008, castro2011} for WENO reconstruction for $r=3$ and $r=4$.

\subsection{Component-wise and characteristic-wise reconstruction}
For system of conservation laws, like the compressible Euler equations, the reconstruction described above can either be done component-wise or characteristic-wise
\cite{zhang2011, qiu2002}. For characteristic-wise decomposition, cell average of the vector of conserved variables $Q$ or the split fluxes $E^{\pm}$ are
transformed into the local characteristic coordinates. The reconstruction is performed on the quantities in the characteristic coordinates and then the reconstructed
quantities are transformed back and used. The transformation to characteristic coordinates and back is done based on the left and right Eigen vectors of the
flux Jacobian $A(Q)$ based on the $Q$ at
the grid point or cell immediately to the left of the point at which the reconstruction is sought, similar to the method labelled U1ZWENO in \cite{zhang2011}.

\subsection{Flux splitting}\label{Sec:FSplittingExp}
Equation \ref{eq:FSplitting} gives the split fluxes($E^{\pm}$) in terms of the flux function ($E$), $Q$ and a parameter $\hat{A}$. 
While calculating approximations to $h^{+}$ and $h^{-}$ at $x_{i+\frac{1}{2}}$, there are different choices for $\hat{A}$.
Choosing $\hat{A} = \alpha$, where $\alpha = \max_{Q}(|\vec{V}| + a)$($a$ is the speed of sound and maximum is taken over all grid points),
leads to the Lax-Freidrich Flux splitting. 

For a less dissipative splitting we choose $\hat{A}$ based on the Roe Flux \cite{roe1981}. Let $Q_L$ and $Q_R$ be left biased and right biased approximations to 
$Q(x_{i+\frac{1}{2}, t_n})$, obtained using WENO interpolation of the same formal order of accuracy ($2r-1$). We refer to \cite{shu2009} for details of WENO Interpolation.
Let $\tilde{Q}$ be the Roe-average state obtained using $Q_L$ and $Q_R$ (see equations 5.41, 5.48, and 5.51 in \cite{laney1998}) and let $A(Q)$ ($=(\partial/\partial Q)E$)
be the flux Jacobian.
Now, we choose $\hat{A} = |A(\tilde{Q})|$. We label this `Roe flux splitting'
and remark that this flux splitting is less dissipative than the Lax-Freidrichs flux splitting.

\subsection{Numerical flux functions}
Equation \ref{eq:FluxFunction} gives the integral of the flux over time in terms of the reconstructed values of $Q$ and a numerical flux function $\hat{E}(l, r)$. Of the
different numerical flux functions available, we present a comparison of results obtained using, the Roe Flux \cite{roe1981}, Roe Flux with Harten Hyman 2 Fix
\cite[p. 266]{harten1983}, Osher's Flux with P-Ordering \cite[Section 12.3.1, p. 393]{toro2013}, Osher's Flux with O-Ordering \cite[Section 12.3.2, p. 397]{toro2013},
AUSM\textsuperscript{+}-up flux \cite{liou2006} and the global Lax-Freidrichs flux.

\subsection{System of equations with viscous fluxes and two-dimensional equations}\label{Sec:extnToViscEqns}
Described previous sections is the procedure for spatial discretisation of hyperbolic conservation law in one space dimension. For equations with viscous fluxes that have
second derivatives of the form
\begin{equation}
  \frac{\partial }{\partial t}Q(x, t) + \frac{\partial}{\partial x}E(Q(x, t)) + \frac{\partial}{\partial x}\Bigg(\mu(x, t) \frac{\partial}{\partial x} \bigg(E_v\big(Q(x, t)\big)\bigg)\Bigg) = 0,
\end{equation} a procedure similar to the one described in sections \ref{Sec:shuOsherConsFinDiffSch} - \ref{Sec:WENOPROC} can be used twice to get the second derivatives. 
The procedure is first applied to calculate the terms 
$\mu(x, t)\partial/\partial x (E_v(Q(x, t)))$. Then the same procedure is again applied on $\mu(x, t)\partial/\partial x (E_v(Q(x, t)))$ to get the second derivative.
Biasing of approximations involving viscous fluxes is not necessary. Therefore, $\hat{A} = 0$ (in equation \ref{eq:FSplitting}) is used for flux splitting. Also, linear
reconstruction ( equation \ref{eq:linHiOrd} - LR) is used for calculating viscous fluxes.

For equations in two space dimensions such as,
\begin{equation}
  \frac{\partial }{\partial t}Q(x, y, t) + \frac{\partial}{\partial x}E(Q(x, y, t)) + \frac{\partial}{\partial y}F(Q(x, y, t)) = 0,
\end{equation}
the same procedure can be used for discretising the $x$ and $y$ derivatives separately.
The resulting semi-discrete form is integrated in time using TVD-RK3 or the Butcher's RK5 method, explained in section \ref{Sec:RKTimeInteg}.

Next we describe the Discontinuous Galerkin Method.

\subsection{Formulation of Discontinuous Galerkin Method}\label{sec:DGMformulation}
The original Discontinuous Galerkin (DG) finite element method was introduced by Reed and Hill \cite{rh} for solving the neutron transport equation which is a linear hyperbolic equation. It was later developed for solving time dependent nonlinear hyperbolic conservation laws as the Runge-Kutta Discontinuous Galerkin (RKDG) method by Cockburn et al. in a series of papers \cite{cs1}, \cite{cs2}, \cite{chs} and \cite{cs3}. The history and development of the DG method is given in the survey paper \cite{cks}.

We now look at solving \eqref{eq:hypConsLaw} using the Discontinuous Galerkin method. We approximate the solution domain by $K$ non overlapping elements whose domain is given by $\mathbf{I}^{k}=[x_{l}^{k},x_{r}^{k}]$. We will approximate the local solution as a polynomial of order $N=N_{p}-1$, where $N_{p}$ is the number of degrees of freedom of the approximation. This is termed to be $\mathbf{P}^{N}$ based Discontinuous Galerkin method. The approximation is given as:

\begin{equation}\label{modalForm}
 Q_{h}^{k}(x,t) = \sum_{n=0}^{N} \hat{Q}^{k}_{n}(t)\psi_{n}^{k}(x) \qquad \forall x\in\mathbf{I}^{k}
\end{equation}

Here, $Q_{h}^{k}(x,t)$ is the approximate local polynomial solution, $\psi_{n}^{k}(x)$ is the local polynomial basis of approximation and $\hat{q}^{k}_{n}(t)$ are the degrees of freedom.

Similarly, we will also approximate the flux $E(Q)$ in the solution domain as given below:

\begin{equation}\label{fluxApprox}
 E_{h}^{k}(Q_{h}^{k}) = \sum_{n=0}^{N} \hat{E}^{k}_{n}(t)\psi_{n}^{k}(x) \qquad \forall x\in\mathbf{I}^{k}
\end{equation}

\noindent We have used the orthonormalized Legendre polynomials as done by Hesthaven et al\cite{hestha1}. The following affine mapping is employed.

\begin{equation}\label{affineMap}
 x(r) = x_{l}^{k} + \frac{1+r}{2}h^{k}, \qquad h^{k} = x_{r}^{k} - x_{l}^{k} \qquad \forall r\in\mathbf{I}=[-1,1]
\end{equation}

\noindent The corresponding recurrence formula for the required orthonormalized Legendre polynomials is given by:

\begin{equation}\label{recFormula}
 r\tilde{P}_{n}(r) = a_{n}\tilde{P}_{n-1}(r) + a_{n+1}\tilde{P}_{n+1}(r),\qquad a_{n} = \sqrt{\frac{n^{2}}{(2n+1)(2n-1)}}
\end{equation}

\noindent with

\begin{displaymath}
 \tilde{P}_{0}(r) = \frac{1}{\sqrt{2}},\qquad \tilde{P}_{1}(r) = \sqrt{\frac{3}{2}}r
\end{displaymath}

\noindent Now the local polynomial basis is given as:

\begin{equation}\label{polyBasis}
 \psi_{n}^{k}(r) = \tilde{P}_{n-1}(r)
\end{equation}

\noindent The degrees of freedom $\hat{Q}^{k}_{n}$ can be advanced in time by the following scheme obtained from the weak form of the governing equation:

\begin{equation}\label{weakFormScheme}
 \frac{d}{dt}\hat{Q}_{h}^{k} = (\mathbf{M}^{k})^{-1}(\mathbf{S}^{k})^{T} \hat{E}_{h}^{k}({Q}_{h}^{k}) - (\mathbf{M}^{k})^{-1} (E^{*}|_{r_{N_{p}}}e_{N_{p}} - E^{*}|_{r_{1}}e_{1})
\end{equation}

\noindent Here, $e_{i}$ is a vector of length $N_{p}$ which has zero entries everywhere except at the $i$th location, and $\mathbf{M}^{k}$ is the local mass matrix which is given as:

\begin{equation}\label{massMatrix}
 \mathbf{M}^{k} = \left[M_{ij}^{k}\right] = \left[\int_{x_{l}^{k}}^{x_{r}^{k}} \psi_{i}^{k}(x) \psi_{j}^{k}(x) \text{dx}\right]
\end{equation}

\noindent and $\mathbf{S}^{k}$ is the local stiffness matrix which is given by:

\begin{equation}\label{stiffnessMatrix}
 \mathbf{S}^{k} = \left[S_{ij}^{k}\right] = \left[\int_{x_{l}^{k}}^{x_{r}^{k}} \psi_{i}^{k}(x) \frac{d\psi_{j}^{k}(x)}{dx} \text{dx}\right]
\end{equation}

\noindent Also, $E^{*}$ is the monotone numerical flux at the interface which is calculated using an exact or approximate Riemann solver. A study of performance of various numerical fluxes for discontinuous Galerkin method has been done in \cite{qks}.
\\
\\
\noindent Now, the semi-discrete scheme given in \eqref{weakFormScheme} is discretized in time by using the TVD Runge-Kutta time discretization introduced in \cite{shu}. We have used a third order TVD Runge-Kutta time discretization for all our calculations. For equations with viscous fluxes of the form

\begin{equation}\label{eqnDGViscous}
  \frac{\partial }{\partial t}Q(x, t) + \frac{\partial}{\partial x}\left(E(Q(x, t))-E_{v}(Q(x, t),\nabla Q(x,t))\right) =  0,
\end{equation}

\noindent we use the local DG (LDG) method as given by Cockburn and Shu in \cite{cs5}. We will solve \eqref{eqnDGViscous} along with

\begin{equation}\label{gradientEqnDG}
 U(x,t) - \nabla Q(x,t) = 0
\end{equation}

\noindent Here $U(x,t)$ can be approximated locally as 

\begin{equation}\label{modalFormGradient}
 U_{h}^{k}(x,t) = \sum_{n=0}^{N} \hat{U}^{k}_{n}(t)\psi_{n}^{k}(x) \qquad \forall x\in\mathbf{I}^{k}
\end{equation}

\noindent Using this, we can obtain the weak form of \eqref{gradientEqnDG} as 

\begin{equation}
 \mathbf{M}^{k}\hat{U}_{m}^{k} = \int_{\mathbf{I}^{k}} Q_{h}(x,t).\nabla\psi_{m}^{k}(x) dx - \int_{\partial \mathbf{I}^{k}} Q_{h}(x,t)\psi_{m}^{k}(x).\hat{n} ds 
\end{equation}

\noindent Each of the above integral is evaluated using an appropriate quadrature rule. The value $Q_{h}(x,t).\hat{n}$ is part of a surface integral and it is taken to be $Q_{h}^{+}.\hat{n}$ where $+$ indicates the discontinuous value outside the element. The value of $\psi_{m}^{k}(x)$ on the surface integral is taken from inside the element. This way, once we obtain $U_{h}(x,t) = \nabla Q(x,t)$, we can find all terms in $E_{v}(Q(x, t),\nabla Q(x,t))=E_{v}(Q(x, t),U(x,t))$. Then \eqref{eqnDGViscous} is written in the weak form similar to \eqref{weakFormScheme} and we can solve the whole system of equations. The numerical flux for $E$ labelled $E^{*}$ in the weak form is obtained using an exact or approximate Riemann solver. The numerical flux for $E_{v}$ labelled $E_{v}^{*}$ is taken to be $E_{v}^{-}$ where $-$ represents the discontinuous value of the solution inside the element. We again use a third order TVD Runge-Kutta time discretization for the solution of the system in time. This completes the LDG formulation.
\\
\\
\noindent Solutions obtained with Discontinuous Galerkin method develop spurious oscillations near discontinuities and a non linear limiter is used to control such oscillations. The common methodology for limiting in Discontinuous Galerkin method is as given below in two steps:
\\
\noindent \textbf{1)} Identify the cells which need to be limited. They are often called troubled cells. \\
\noindent \textbf{2)} Replace the solution polynomial in the troubled cell with a new polynomial that is less oscillatory but with the same cell average and order of accuracy.
\\
\\
\noindent For the first step, we have used the KXRCF troubled cell indicator for all the calculations done in this paper as it is rated highly by Qiu and Shu in \cite{qs2} on the basis of it's performance in detecting the discontinuities in various test problems. The second step is where we do the limiting process. We have used the so called simple WENO limiter developed by Zhong and Shu \cite{zs} for all the calculations done in this paper.

\subsection{Labelling the numerical methods}
We will denote component-wise reconstruction methods using weights given by equations (\ref{eq:linHiOrd})
, (\ref{eq:omegaWenoZ}) as LR,
ZW respectively and characteristic-wise reconstruction as LCDLR,
LCDZW, respectively.  The number
following these labels is used to indicate the formal order of accuracy of the reconstruction. 
We use a prefix FD and FV to denote conservative finite difference and finite volume methods respectively
Similarly we denote the Discontinuous Galerkin methods using DG $P^n$ label,
where the $n$ indicates the degree of the basis polynomial.
To indicate the flux splitting or flux function used, we add one of the
following suffixes:
\begin{itemize}
    \item `-ROE' for Roe Flux
    \item `-ROEHH2' for ROE Flux with Harten Hyman 2 entropy fix,
    \item `-LF' for Lax-Freidrichs flux function,
    \item `-AUSM' for the AUSM\textsuperscript{+}-up flux,
    \item `-OshP' for Osher's P-Ordering flux.
    \item `-OshO' for Osher's O-Ordering flux.
    \item `-C' suffix indicates central scheme was used.
    \end{itemize}
Therefore, FDLR7-C indicates linear reconstruction with formal order of accuracy of 7 with $\hat{A}=0$ is used. 
FVLCDZW3-Roe indicates that finite volume method with characteristic-wise WENO-NP3 reconstruction and Roe flux is
used while FDZW5-LF indicates that finite difference method with component-wise ZWENO5 reconstruction with Lax-Freidrichs flux splitting is used.

\subsection{Runge Kutta time discretisation}\label{Sec:RKTimeInteg}
Consider the equation
\begin{equation}
  \label{eq:timInt}
  \frac{d}{dt} Q = L(Q).
\end{equation}
The simple forward Euler time discretisation between two time levels $t_n$ and $t_{n+1}$ separated by $\Delta t$ is given by
\begin{equation}
  \label{eq:fowEulDisc}
  Q^{n+1} = Q^{n} + \Delta t L(Q^{n}).
\end{equation}
A three stage third order TVD (Total Variation Diminishing) or SSP (Strong Stability Preserving) \cite{gottlieb2001} Runge-Kutta discretisation is given by
\begin{align}
  Q^{(1)} = Q^{n} + \Delta t L(Q^{n}),~~&
  Q^{(2)} = \frac{3}{4} Q^n + \frac{1}{4}Q^{(1)} + \frac{1}{4}\Delta t L(Q^{(1)}),\\
  Q^{n+1} = \frac{1}{3} Q^n &+ \frac{2}{3}Q^{(1)} + \frac{2}{3}\Delta t L(Q^{(2)}).
\end{align}
The Runge Kutta discretisation described above are used to advance in time from $t_n$ to $t_{n+1}$.
\subsection{Refinement near shock using overset mesh}\label{Sec:OversetDesc}
To improve solution accuracy near stationary shocks, we use a finer overset mesh. Next, we briefly describe the procedure employed for obtaining numerical solutions using overset mesh.
\subsubsection{Conservative coupling procedure for the finite difference scheme}
\begin{figure}[!htbp]
\begin{center}
  \includegraphics[width=0.65\textwidth]{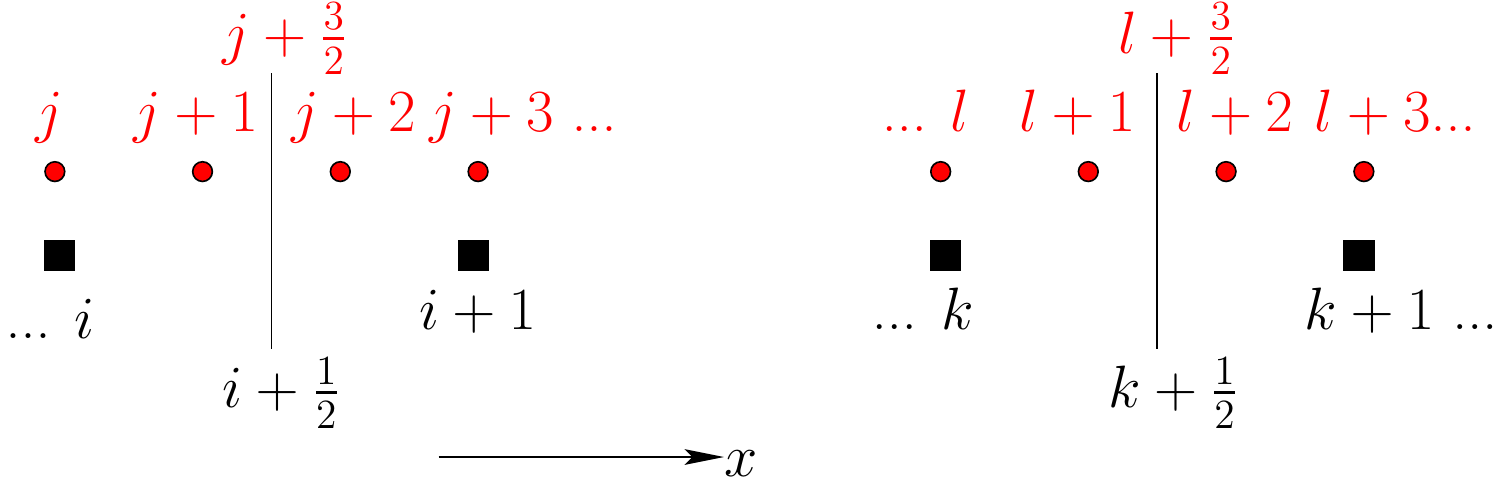}
  \caption{ Overset mesh with coarse (black, square grid points) and fine components (red, circular grid points), with left ($i+\frac{1}{2}$) and right ($k+\frac{1}{2}$) coupling interfaces shown.}
  \label{fig:overset}
\end{center}
\end{figure}
As shown in figure \ref{fig:overset}, for the finite difference method, we use an overset mesh with coarse (black, square grid points) and fine components (red, circular
grid points). The grid point spacing (GPS or $\delta x$) of the finer mesh component is chosen so as to have coupling interfaces like $i+\frac{1}{2}$, $j+\frac{3}{2}$
and overlapping mesh points, like $i$, $j$, $i+1$, $j+3$. On the coarse mesh component, grid points from $i+1$ to $k$ are fringe points and the remaining grid points are
discretisation points, except for the ghost points used for boundary condition application. On the finer mesh component, grid points $j+2$ to $l+1$ are the 
discretisation points and the remaining are fringe points. State in the fringe points in the coarse mesh (like $i+1$) is copied directly from the corresponding
overlapping grid points ($j+3$) in the fine mesh. State in the fringe points in the fine mesh (like $j$, $j+1$) is obtained using an interpolation polynomial,
based on the data from nearest neighbouring grid points in the coarse mesh component. The degree of the polynomial used for interpolation is equal to the formal order of accuracy (F.O.A) of the finite
difference scheme used.

To have a conservative coupling between the coarse and the finer mesh components, a unique numerical flux $\hat{h}^{\pm}$ (see equation \ref{eq:fluxReconsApprox})
must be used at the left ($i+\frac{1}{2}$,$j+\frac{3}{2}$) and right ($l+\frac{1}{2}$, $l+\frac{3}{2}$) coupling interfaces \cite{cheng2013, cheng2016}. This numerical
flux can be calculated using either the coarse mesh component or fine mesh component or any convex combination of them. Using numerical flux from the fine component at the
left coupling interface and the numerical flux from the coarse mesh component at the right coupling interface seems to give good results. 
\subsubsection{Overset mesh method for the DG scheme}
\begin{figure}[!htbp]
\begin{center}
  \includegraphics[width=0.35\textwidth]{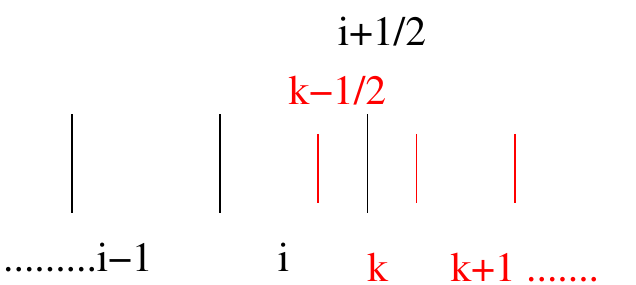}
  \caption{ Overset mesh with coarse (black) and fine components (red).}
  \label{fig:oversetDG}
\end{center}
\end{figure}

A typical one-dimensional overset mesh for DGM is shown in Figure \ref{fig:oversetDG}. Here the element $i$ in the coarse mesh (black) and the element $k$ in the fine mesh (red) are overlapping. While advancing the solution in time in the coarse mesh, we find the solution at $i+1/2$ in the fine mesh (by locating it appropriately in the local coordinate system of the fine mesh) and apply it as the boundary condition for calculating the numerical flux at $i+1/2$. This procedure is followed as given by Galbraith et al\cite{gbot} for two-dimensional meshes. Similarly for advancing the solution in the fine mesh, we apply the coarse mesh solution at $k-1/2$ as the boundary condition. This is again used to calculate the numerical flux at $k-1/2$. This procedure gives good results for using overset meshes with DGM.

\section{Verification}\label{Sec:Verification}
We verify the implementation of the WENO and DG schemes using the problems of the burgers equations with source term and the isentropic Euler vortex problem \cite{spiegel2015}.
\subsection{Burgers equation with source term}
We solve
\begin{equation}
  \frac{\partial}{\partial t} u(x, t) + \frac{\partial}{\partial x}\bigg(\frac{1}{2}u^2(x, t)\bigg) = -8x^7u(x, t), 0\leq x \leq 1
\end{equation}
with the boundary conditions $u(0, t) = 2.0$, $u(1, t) = 1.0$ and the initial conditions $u(x, 0) = 2.0-x$ to steady state. The steady state solution is $u_s(x) = 2.0 - x^8$.
The flux function, $E(x, t) = u^2(x, t)/2$ and the flux splitting is $E^{+} = u^2(x, t)/2$ and $E^{-} =0$ (LB or left biased splitting). The FDZW7-LB
and DG $P^7$-LB scheme with TVD-RK3 time discretisation was used to obtain the numerical solutions for different grid point spacings(GPS = $\Delta x$) or different 
cell sizes. Table 
\ref{tab:burgersL1Error} has the L1 errors and the
observed order of accuracy.
\begin{table}[!htbp]
  \begin{center}
  \caption{$L_1$ errors of numerical solutions obtained using FDZW7-LB and DG schemes for different GPS/cell size and observed order of accuracy.}
\label{tab:burgersL1Error}
  \begin{tabular}{|>{\centering\arraybackslash}m{0.1\textwidth}|>{\centering\arraybackslash}m{0.14\textwidth}|>{\centering\arraybackslash}m{0.08\textwidth}|>{\centering\arraybackslash}m{0.14\textwidth}|>{\centering\arraybackslash}m{0.08\textwidth}|}
      \hline
      \multirow{2}{5em}{GPS/Cell size} & \multicolumn{2}{c|}{FDZW7}&\multicolumn{2}{c|}{DG $P^7$}\\\cline{2-5}
                                       &$L_1\text{ error}$ $\times 10^{-11}$ & order & $L_{1}\text{ error}$ $\times 10^{-11}$ & order \\ \hline 
      1/25 & 120870.8 & - & 23426.1 & - \\ \hline
      1/50 & 851.6 & 7.14 & 162.35 & 7.17 \\ \hline
      1/75 & 49.2 & 7.03 & 9.2153 & 7.075 \\ \hline
      1/100 & 7.0 & 6.77 & 1.2834 & 6.853 \\ \hline
    \end{tabular}
  \end{center}
\end{table}
\subsection{Isentropic Euler Vortex Problem \cite{spiegel2015}}
We solve the two-dimensional Euler equations, which are
\begin{equation}
  \label{eq:Euler2dDiff}
  \frac{\partial Q}{\partial t} + \frac{\partial E}{\partial x} + \frac{\partial F}{\partial y} = 0
\end{equation}
where
\begin{equation}
  \label{eq:Euler2dDiffTerms}
~Q = \begin{bmatrix}\rho \\ \rho u \\ \rho v \\ \rho e_t\end{bmatrix}, ~E =  \begin{bmatrix}\rho u \\ \rho u^2 + p \\ \rho uv \\ (\rho e_t + p)u\end{bmatrix},
~F =  \begin{bmatrix}\rho v \\ \rho vu \\ \rho v^2 +p \\ (\rho e_t + p)v\end{bmatrix},
 ~ e_t = \frac{p}{\rho(\gamma -1)} + \frac{1}{2}\left(u^2 + v^2 \right)
\end{equation}
The initial conditions are an isentropic vortex perturbation added to a uniform flow in the positive $x$ direction and is given by:
\begin{align}
  u(x, y, 0) &= u_0 - 5 e^{(1-r^{2})} \frac{y-y_{0}}{2\pi},\\
  v(x, y, 0) &= 5 e^{(1-r^{2})} \frac{x-x_{0}}{2\pi},\\
  \rho(x, y, 0) &= \left(1 - \left(\frac{\gamma - 1}{16\gamma \pi^{2}}\right) 25 e^{2(1-r^{2})}\right)^{\frac{1}{\gamma-1}},
\end{align}
with $p(x, y, 0) = (\rho(x, y, 0))^{\gamma}$ and $r=\sqrt{(x-x_{0})^{2}+(y-y_{0})^{2}}$. The parameter values chosen are $x_{0}=7.0$, $y_{0}=0$, $\beta=2.0$, $u_0 = 1.0$
and $\gamma = 1.4$. The computational domain is a square of dimensions $14 \text{units} \times 14 \text{units}$ with $0 \leq x \leq 14$ and $-7 \leq y \leq 7$. Periodic
boundary conditions are applied along the $x$ and $y$ directions. The FDZW5-LF and DG $P^4$-LF method with Lax-Friedrichs flux splitting are used to obtain numerical solution 
at $t = 14.0$ units (one time period).
For time discretisation Butchers six stage and fifth order RK scheme with time step ${\Delta t= 0.07 \times \Delta x}$ was used with the FDZW5-LF scheme. TVD-RK3
scheme with $\Delta t = 0.1 \times (\Delta x)^{5/3}$ for the DG $P^4$ scheme. This problem was run for meshes with grid point
spacings ($GPS=\Delta x = \Delta y$) of $1/25, 1/50, 1/75, 1/100, 1/150, 1/175, 1/200, \text{ and } 1/225$. The $L_1$
errors for meshes with different GPS and the observed order of accuracy are given in table~\ref{tab:eulerIsenVortErrs}.
\begin{table}[!htbp]
  \begin{center}
    \caption{$L_1$ errors of total energy density ($\rho e_t$) obtained using FDZW5-LF and DG schemes for different GPS/cell size and observed order of accuracy.}
\label{tab:eulerIsenVortErrs}
  \begin{tabular}{|>{\centering\arraybackslash}m{0.1\textwidth}|>{\centering\arraybackslash}m{0.14\textwidth}|>{\centering\arraybackslash}m{0.08\textwidth}|>{\centering\arraybackslash}m{0.14\textwidth}|>{\centering\arraybackslash}m{0.08\textwidth}|}
      \hline
      \multirow{2}{5em}{GPS/Cell size} & \multicolumn{2}{c|}{FDZW5-LF}&\multicolumn{2}{c|}{DG $P^4$-LF}\\\cline{2-5}
                                       &$L_1\text{ error}$ $\times 10^{-11}$ & order & $L_{1}\text{ error}$ $\times 10^{-11}$ & order \\ \hline 
      1/25 & 235635.2 & - & 15415.2 & - \\ \hline
      1/50 & 6746.7 & 5.13 & 421.64 & 5.19 \\ \hline
      1/75 & 876.3 & 5.03 & 52.345 & 5.145 \\ \hline
      1/100 & 203.2 & 5.08 & 12.346 & 5.02 \\ \hline
      1/150 & 26.7 & 5.00 & 1.6432 & 4.97 \\ \hline
      1/200 & 6.4 & 4.96 & 0.4124 & 4.805 \\ \hline
    \end{tabular}
  \end{center}
\end{table}

\section{Post shock oscillations and mass conservation error}\label{Sec:postShockOscillations}
When there are discontinuities in the solution, the shock capturing methods described above can have issues like mass conservation error, post shock oscillations,
convergence stalling. In this paper, we focus on the issues of the post shock oscillation and mass conservation error. To demonstrate these issues, we
use the problem of a moving normal shock modelled using one-dimensional Euler equations.
\subsection{One-dimensional Euler equations}
Consider the system of equations
\begin{equation}
  \label{eq:oneDimEulerEqn}
  \frac{\partial}{\partial t}Q(x, t) + \frac{\partial}{\partial x}E(x, t) = 0,
\end{equation}
where $Q = [\rho, \rho u , \rho e_t]^T$, $E(x, t) = [\rho u, \rho u^2, (\rho e_t + p)u]^T$, with $\gamma = 1.4$,

We solve equations (\ref{eq:oneDimEulerEqn}), for $0 \leq x \leq L$ with initial conditions
\begin{equation}
  \label{eq:oneDNormShk}
  Q(x, 0) = \begin{cases}
    Q_{BS} & x<x_S \\
    Q_{AS} & x\geq x_S\\
  \end{cases},
\end{equation}
where $(\rho_{BS}, u_{BS}, p_{BS}) = (\gamma, M+u_S, 1.0)$,
$$(\rho_{AS}, u_{AS}, p_{AS}) = \left(\frac{(\gamma+1)M^2 \rho_{BS}}{(\gamma-1)M^2 + 2}, \frac{\rho_{BS}u_{BS}}{\rho_{AS}}+u_S, \frac{p_{BS}(2\gamma M^2 - (\gamma -1))}{\gamma+1}
\right),$$
with supersonic inflow conditions at $x=0$ and subsonic outflow conditions with a back pressure ($p_{back}$) equal to $p_{AS}$ at $x=L$ as boundary conditions. 
These conditions correspond to a normal shock moving with a velocity of $u_s$. We choose a mesh with cell size ($ = \Delta x$) of $1/100$, with number of cells equal to
$L/(\Delta x)$.
Next, we demonstrate the post shock oscillations and mass conservation error using the first order numerical solutions, using different numerical flux functions and
values for parameters $M, u_S, x_S$.


\subsection{Error in numerical solutions obtained using first order schemes}
\begin{figure}[!htbp]
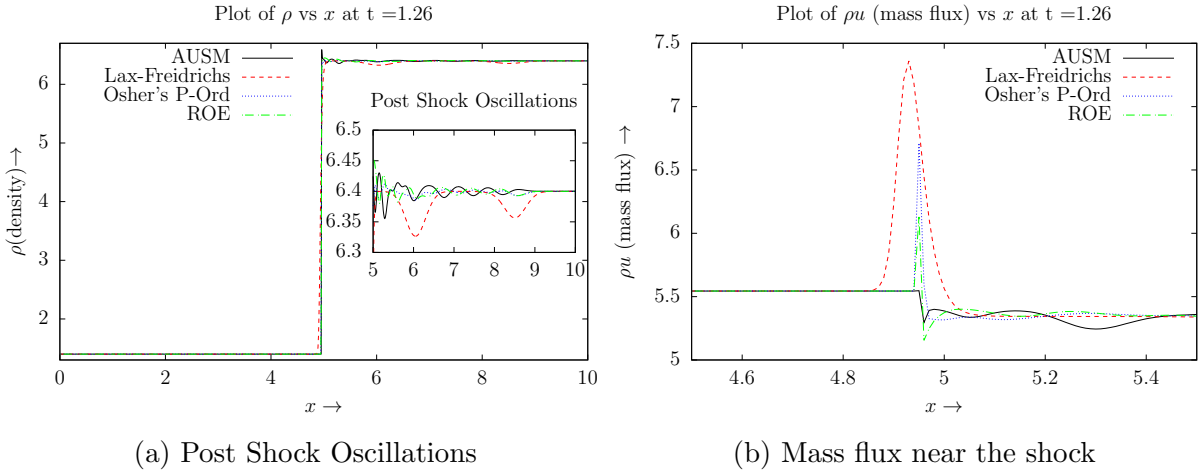

  \centering
  \begin{subfigure}{.5\textwidth}
    \centering
  \scalebox{.65}{\input{images/StM4SpM0_04DensityPSOsc.imgtex}}
  \caption{Post Shock Oscillations}
  \label{fig:StM4SpM0_04PSOsc}
\end{subfigure}%
\begin{subfigure}{.5\textwidth}
  \centering
  \scalebox{.65}{\input{images/StM4SpM0_04DensityMFE.imgtex}}
  \caption{Mass flux near the shock}
  \label{fig:StM4SpM0_04MFE}
\end{subfigure}
\caption{Post shock oscillations and Mass conservation error in numerical solutions (F.O.A = 1) for $M=4$, $u_s = 0.04$ using FVLR1 (finite volume, linear reconstruction ,
first order)}
\label{fig:StM4SpM0_04PSOAndMFE}
\end{figure}

For demonstration, we choose the parameter values $M = 4$, $u_S = -0.04, L = 10$ and $x_S = 5$. Figure \ref{fig:StM4SpM0_04PSOsc}
has plots of $\rho$ vs
x at $t = 1.26$ units for numerical solutions obtained using AUSM\textsuperscript{+}-up, ROE, P-Ordering Osher's and the global Lax-Freidrichs flux functions using
FVLR1
scheme,
which show the post shock oscillations. Figure \ref{fig:StM4SpM0_04MFE}
has plots of $\rho u$ vs x at $t = 1$ units, obtained
using different numerical flux functions. These figures show the momentum density or mass flux spike (non monotonic variation in the mass flux). Both the post shock
oscillations and the mass flux spike are artefacts of the numerical solutions.
We seek to quantify these errors. 

We know that across the moving shock ($0 \leq x \leq L$), $CAS_1, CAS_2, CAS_3$ defined by
\begin{equation}
  \label{eq:movShkInvariants}
  CAS_1 = \rho(u - u_S), CAS_2 = \rho(u-u_S)^2 + p, CAS_3 = \frac{p}{\gamma-1} + \frac{1}{2}\rho(u-u_S)^2,
\end{equation}
are constant. 
Based on the invariant $CAS_1$, we define the total mass conservation error percentage (CEP) in the numerical solution at a time $t_n$ as
\begin{equation}
  \label{eq:totalMCEDef}
  \text{CEP}_n = \int\limits_{x=0}^{x=L}\frac{\rho(x, t_n)(u(x, t_n)-u_S) - \rho_{BS}u_{BS}}{\rho_{BS}u_{BS}}dx \times 100.
\end{equation}
Equation \ref{eq:totalMCEDef} can be written in terms of cell averages as
\begin{equation}
  \label{eq:totalMCEDisc}
  \text{CEP}(t_n) = \text{CEP}_n = \sum\limits_{i=0}^{\text{num. cells}}\frac{\overline{\rho u}_i^n - \bar{\rho}_i^nu_S - \rho_{BS}u_{BS}}{\rho_{BS}u_{BS}}\Delta x \times 100
\end{equation}
\begin{figure}[!htbp]
  \centering
  \scalebox{.9}{\input{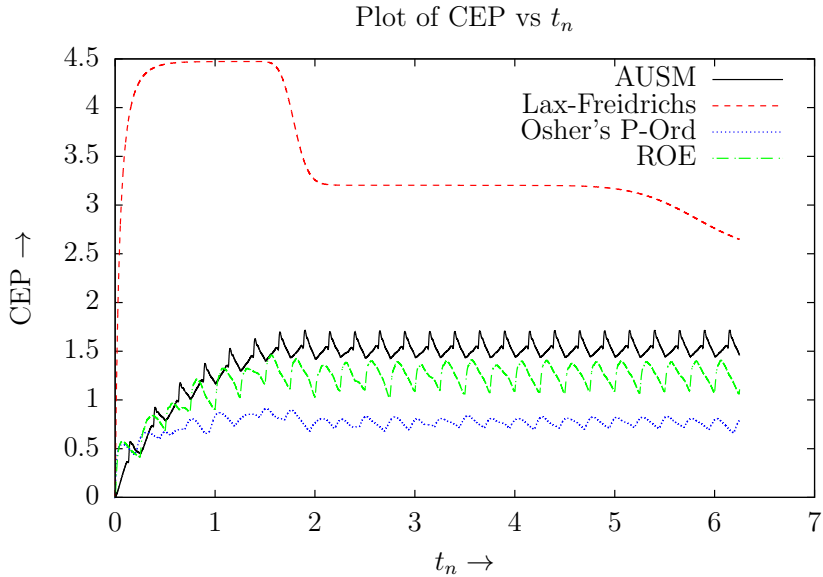}}
  \caption{Total mass conservation error percentage (CEP) for $M=4$, $u_s = -0.04$ using FVLR1}
\label{fig:StM4SpM0_04CEP}
\end{figure}
\begin{figure}[!htbp]
  \centering
  \scalebox{.9}{\input{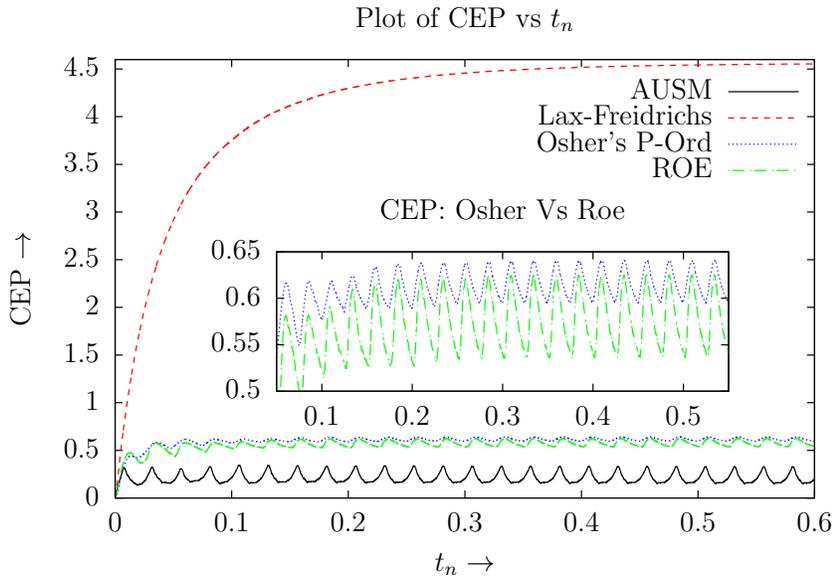}}
  \caption{Total mass conservation error percentage (CEP) for $M=4$, $u_s = -0.4$ using FVLR1}
\label{fig:StM4SpM0_4CEP}
\end{figure}

Figure \ref{fig:StM4SpM0_04CEP}
has plots of CEP$(t_n)$ vs $t_n$ ($n$ = time step number) in numerical solutions obtained using different numerical flux functions. We can see that the solution obtained
using Osher's P-Ordering flux has less error than the that
of ROE flux which is consistent with what is reported in literature \cite{roberts1990}. This however changes with changes in $u_s$, for instance for $u_s = -0.4$, shown in
figure \ref{fig:StM4SpM0_4CEP}
For $u_s = -0.4$, the order of error is reversed, with least error obtained using AUSM\textsuperscript{+}, followed by ROE flux
and Osher's P-Ordering flux. In both cases ($u_s = -0.04$  and $u_s = -0.4$), using the Lax-Freidrichs flux leads to the highest error.

Next, the errors in numerical solutions obtained using high-order schemes for the same problem is discussed.

%

\subsection{Error in numerical solutions obtained using high order schemes}
\begin{figure}[!htbp]
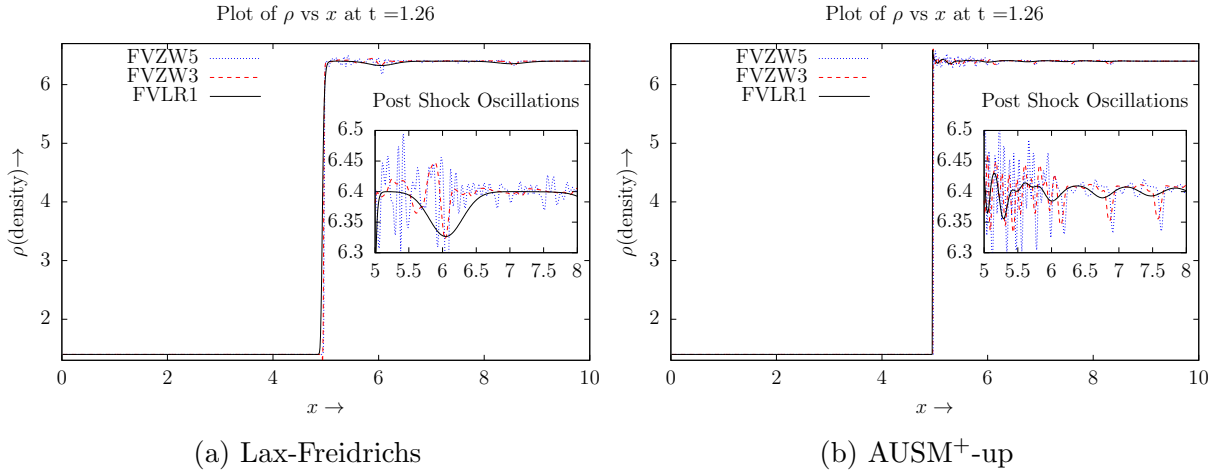

  \centering
  \begin{subfigure}{.5\textwidth}
    \centering
  \scalebox{.65}{\input{images/StM4SpM0_04DensityPSOscLaxFrCmp.imgtex}}
  \caption{Lax-Freidrichs}
  \label{fig:StM4SpM0_04PSOscLaxFr}
\end{subfigure}%
\begin{subfigure}{.5\textwidth}
  \centering
  \scalebox{.65}{\input{images/StM4SpM0_04DensityPSOscAUSMPlusCmp.imgtex}}
  \caption{AUSM\textsuperscript{+}-up}
  \label{fig:StM4SpM0_04PSOscAUSMPl}
\end{subfigure}
\caption{Comparison of post shock oscillations obtained using schemes with different F.O.A for AUSM\textsuperscript{+}-up splitting and Lax-Freidrichs flux
for $M=4$, $u_s = -0.04$ using FVZW and FVLR}
\label{fig:StM4SpM0_04PSOscCmp}
\end{figure}
\begin{figure}[!htbp]
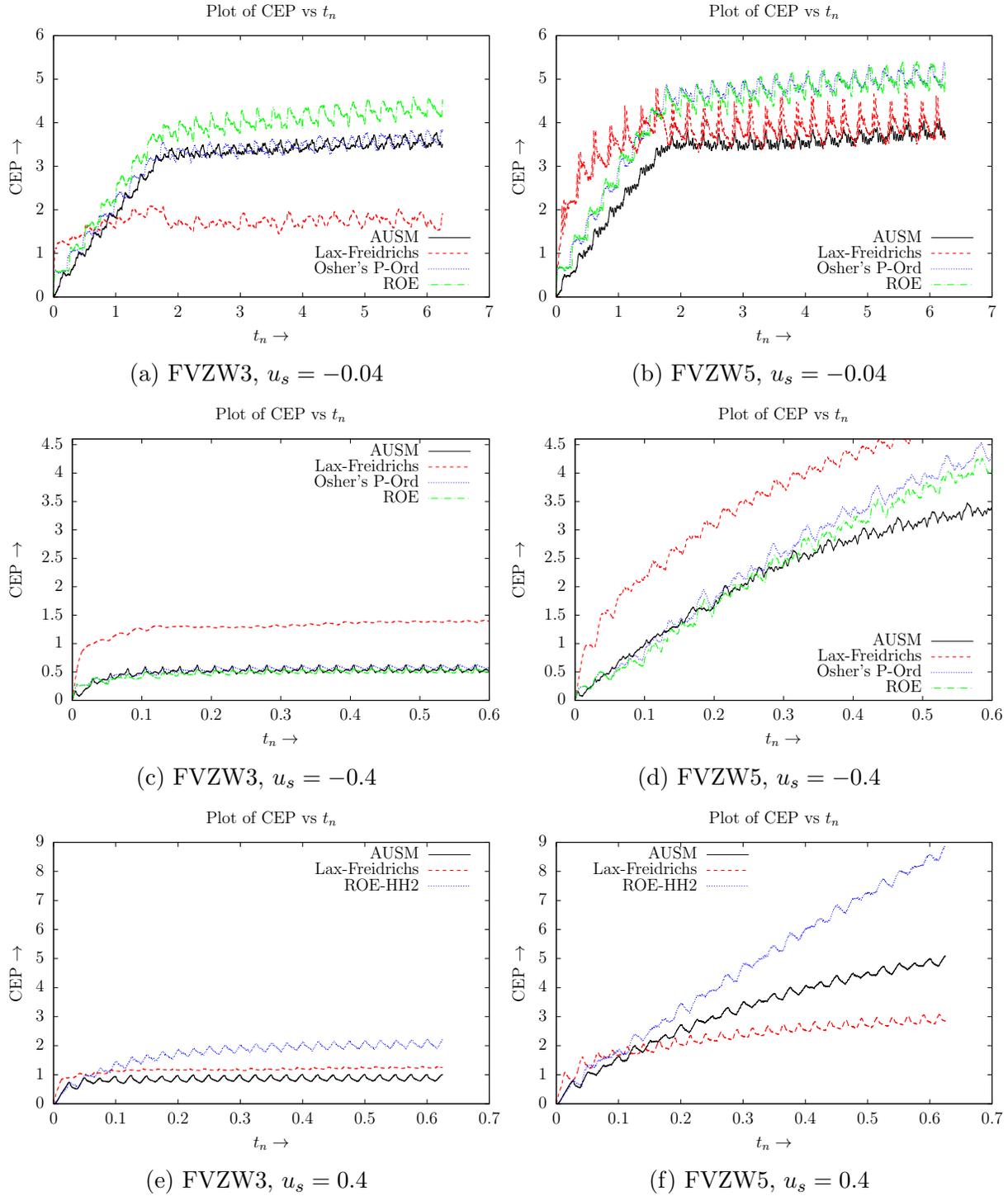

  \centering
  \begin{subfigure}{.5\textwidth}
    \centering
  \scalebox{.65}{\input{images/StM4SpM0_04TotalMFEW3.imgtex}}
  \caption{FVZW3, $u_s = -0.04$}
  \label{fig:StM4SpM0_04TotalMFEW3}
\end{subfigure}%
\begin{subfigure}{.5\textwidth}
  \centering
  \scalebox{.65}{\input{images/StM4SpM0_04TotalMFEW5.imgtex}}
  \caption{FVZW5, $u_s = -0.04$}
  \label{fig:StM4SpM0_04TotalMFEW5}
\end{subfigure}
  \centering
  \begin{subfigure}{.5\textwidth}
    \centering
  \scalebox{.65}{\input{images/StM4SpM0_4TotalMFEW3.imgtex}}
  \caption{FVZW3, $u_s = -0.4$}
  \label{fig:StM4SpM0_4TotalMFEW3}
\end{subfigure}%
\begin{subfigure}{.5\textwidth}
  \centering
  \scalebox{.65}{\input{images/StM4SpM0_4TotalMFEW5.imgtex}}
  \caption{FVZW5, $u_s = -0.4$}
  \label{fig:StM4SpM0_4TotalMFEW5}
\end{subfigure}
  \centering
  \begin{subfigure}{.5\textwidth}
    \centering
  \scalebox{.65}{\input{images/StM4SpMPl0_4TotalMFEW3.imgtex}}
  \caption{FVZW3, $u_s = 0.4$}
  \label{fig:StM4SpMPl0_4TotalMFEW3}
\end{subfigure}%
\begin{subfigure}{.5\textwidth}
  \centering
  \scalebox{.65}{\input{images/StM4SpMPl0_4TotalMFEW5.imgtex}}
  \caption{FVZW5, $u_s = 0.4$}
  \label{fig:StM4SpMPl0_4TotalMFEW5}
\end{subfigure}
\caption{Comparison of CEP in numerical solutions obtained using FVZW3, FVZW5 methods with different numerical flux functions, for $u_s = -0.4, -0.04$ and $0.4$ }
\label{fig:StM4TotalMFEW3W5Cmp}
\end{figure}
Figure \ref{fig:StM4SpM0_04PSOscCmp}
has plots of density ($\rho$) vs $x$ at $t=1.26$, obtained using FVLR1, FVZW3 FVZW5 schemes and  Lax-Freidrichs, AUSM\textsuperscript{+}-up
fluxes, for $u_s = -0.04$ . As the formal order of accuracy of the scheme increases, both the amplitude and wave number of the oscillations increase.
A similar trend is observed for the ROE flux and Osher's P-Ordering flux. Figure \ref{fig:StM4TotalMFEW3W5Cmp} has plots of CEP, for different
values of $u_s$, different numerical flux functions, for FVZW3 and FVZW5 methods.
From the plots in figures \ref{fig:StM4SpM0_04CEP},
\ref{fig:StM4SpM0_4CEP}, and
\ref{fig:StM4TotalMFEW3W5Cmp},
it seems that no one numerical flux function consistently leads to less error. For instance in the case
of $u_s=0.4$, for the FVZW5 method, Lax-Freidrichs flux function leads to error (see figure \ref{fig:StM4SpMPl0_4TotalMFEW5}), less than that due to the
AUSM\textsuperscript{+}-up and the ROE flux with Harten-Hyman 2 \cite[p. 266]{harten1983} entropy fix.

Additionally, for $u_s = 0.4$, the ROE and Osher's P-Ordering fluxes fail to produce numerical solutions when FVZWE3 or FVZW5 methods are used as `NAN' is
produced in the course of computations. Using the Harten-Hyman 2 \cite[p. 266]{harten1983} entropy fix for the ROE flux, rectified the problem and a numerical solution
was obtained, but with an error higher than that of Lax-Freidrichs and AUSM\textsuperscript{+}-up fluxes (see Figures \ref{fig:StM4SpMPl0_4TotalMFEW3},
\ref{fig:StM4SpMPl0_4TotalMFEW5}).

The problem of encountering `NAN' in using ROE and Osher's P-Ordering fluxes can be tackled by doing a characteristic-wise reconstruction. 


\subsection{The importance of characteristic decomposition}
\begin{figure}[!htbp]
  \centering
  \begin{subfigure}{.5\textwidth}
    \centering
  \scalebox{.65}{\input{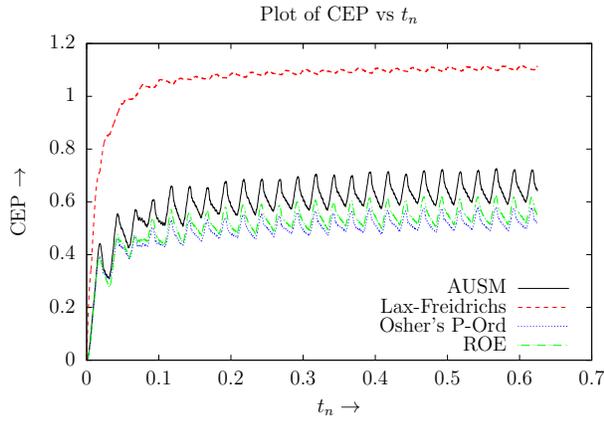}}
  \caption{LCDZW3}
  \label{fig:StM4SpMPl0_4TotalMFELCDW3}
\end{subfigure}%
\begin{subfigure}{.5\textwidth}
  \centering
  \scalebox{.65}{\input{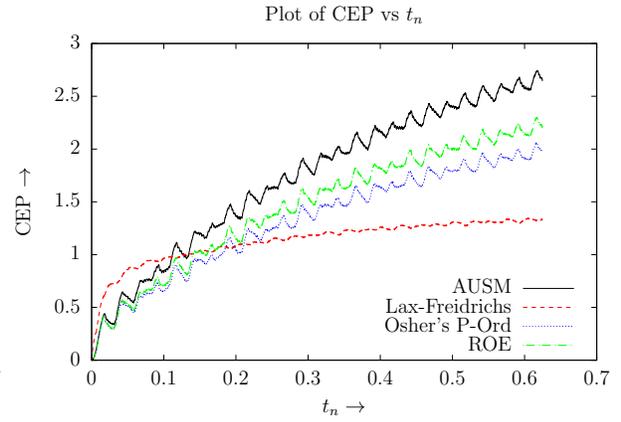}}
  \caption{LCDZW5}
  \label{fig:StM4SpMPl0_4TotalMFELCDW5}
\end{subfigure}
  \centering
  \begin{subfigure}{.5\textwidth}
    \centering
  \scalebox{.65}{\input{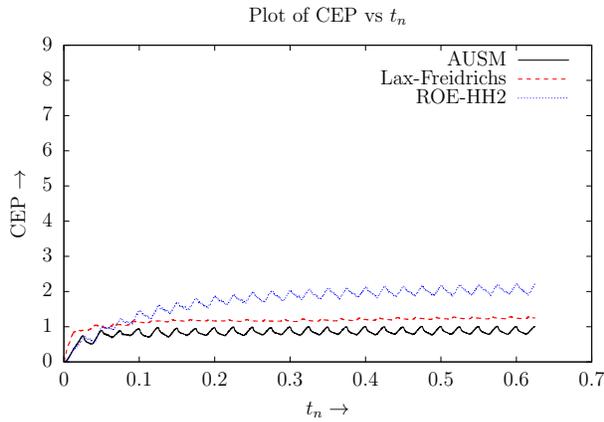}}
  \caption{DG - P2, $u_s = 0.4$}
  \label{fig:StM4SpMPl0_4TotalMFEW3DG}
\end{subfigure}%
\begin{subfigure}{.5\textwidth}
  \centering
  \scalebox{.65}{\input{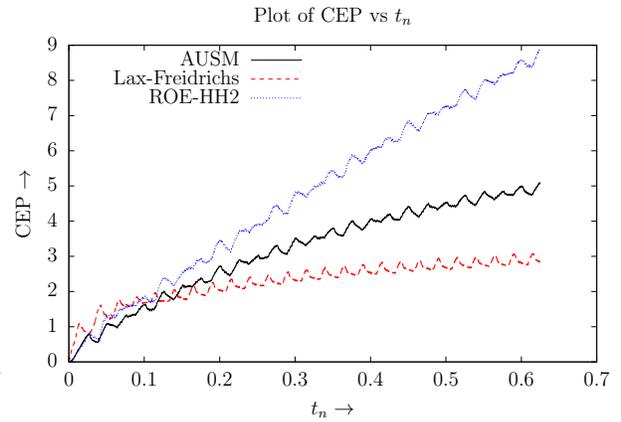}}
  \caption{DG - P4, $u_s = 0.4$}
  \label{fig:StM4SpMPl0_4TotalMFEW5DG}
\end{subfigure}
\caption{Comparison of CEP in numerical solutions obtained using LCDZW3, LCDZW5, DG $P^2$, DG $P^4$  and different numerical flux functions, for $u_s = 0.4$.
}
\label{fig:StM4SpMPl0_4TotalMFELCDW}
\end{figure}
As mentioned before, doing a characteristic-wise reconstruction \cite{qiu2002, peng2019, zhang2011}, will produce better and less oscillatory results when compared to doing
component wise reconstruction. An example illustrating this point in the extreme is the problem discussed above, for $u_s = 0.4$. Doing a component wise ZWENO3 or
ZWENO5 reconstruction, and using the ROE or Osher's P-Ordering flux will lead to `NAN's in the computations. This can be rectified by doing a characteristic-wise
reconstruction. The characteristic-wise reconstruction also reduces the error or `CEP' as can be seen in figures \ref{fig:StM4SpMPl0_4TotalMFEW3},
\ref{fig:StM4SpMPl0_4TotalMFEW5} and \ref{fig:StM4SpMPl0_4TotalMFELCDW}. Errors similar to FDLCDZW are obtained using DG with simple WENO limiter (with 
characteristic-wise limiting), as shown in figures \ref{fig:StM4SpMPl0_4TotalMFEW3DG} and \ref{fig:StM4SpMPl0_4TotalMFEW5DG}.

In summary, the order of performance of numerical flux function seems problem dependent and while using high order reconstruction, doing a characteristic wise
reconstruction is critical to get results. As the formal order of accuracy increases, results produced using Lax-Freidrichs flux seem to become better than that obtained using high resolution fluxes, as is
evident from figure \ref{fig:StM4SpMPl0_4TotalMFELCDW}. Next, we study the mass conservation error in steady state
numerical solutions having shocks.

\section{Mass conservation error in numerical solutions with stationary shocks}\label{Sec:massFluxErr}
We solve the Euler equations \ref{eq:oneDimEulerEqn} in $0 \leq x \leq 1$ with the initial conditions given in equation (\ref{eq:oneDNormShk}), with $u_s = 0$.
We obtain the numerical solutions using TVD-RK3 time discretisation, using FDLCD and DG methods.  Evaluating the spatial derivative
($(\partial E)/(\partial x)$) using the ROE flux leads to a zero value for the spatial derivative whereas for the Lax-Freidrichs flux it leads to a non-zero value.
The also leads to an error in
the mass flux($\rho u$) near the shock. Across a normal shock, $\rho u$ should be constant, where as in the numerical solution obtained using Lax-Freidrichs splitting 
$\rho u$ is not constant. We define the percentage error of a conserved variable($q$) at grid point $i$ and time $t_n$ as
\begin{equation}
  \label{eq:percentageError}
  PE(q_i^n) = \frac{|q(0, 0) - q_i^n|}{q(0, 0)} \times 100.
\end{equation}
Here $q_i^n$ is the value of the conserved variable $q$ at grid point $x_i$, at time $t_n$ and $q(0, 0)$ is the value of $q(x,t)$ at $x=0,t=0$. For a steady state solution of the one-dimensional Euler equations, $q$ can be one of $\rho u, (\rho u^2 + p)$ and
$(\rho e_t + p)u$ as they are conserved, whereas in the numerical
solution they
are not conserved.

\begin{figure}[!htbp]
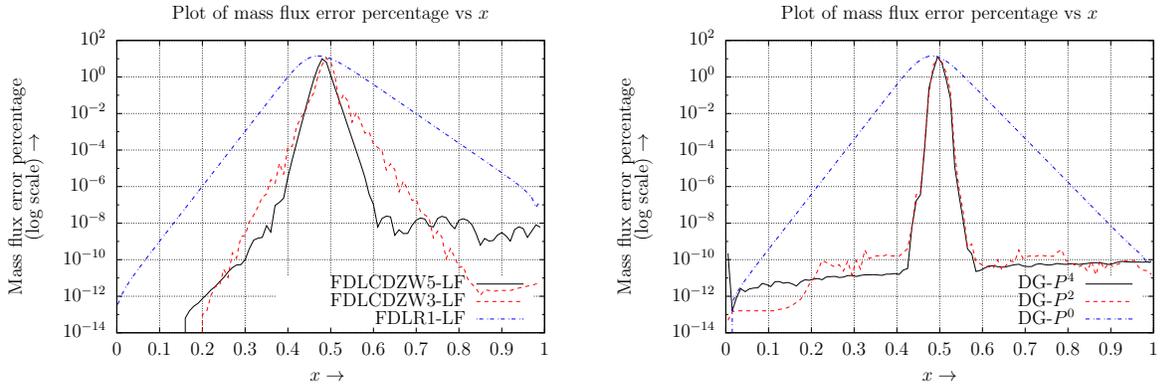

  \centering
  \begin{subfigure}{.5\textwidth}
    \centering
  \scalebox{.6}{\input{images/massFluxErr/mach2MassErrWENODiffFOA.imgtex}}
  \caption{FDLCDZW5,3 and FDLR1 with LF splitting}
  \label{fig:massFluxE}
\end{subfigure}%
\begin{subfigure}{.5\textwidth}
  \centering
  \scalebox{.6}{\input{mach2MassErrDGDiffFOA.imgtex}}
  \caption{DG - $P^4, P^2, P^0$ with LF splitting}
  \label{fig:massFluxErrM2DGP4}
\end{subfigure}
\caption{Mass flux error percentages across a Mach 2 shock at $t=100.0$, for schemes with different formal order of accuracies.}
\label{fig:massErrorDiffFOA}
\end{figure}

Figure \ref{fig:massErrorDiffFOA} has plots of $PE(\rho u_i^n)$ vs $x_i$ at $t_n=100.0$ (mass flux error percentages), for schemes with different formal order of accuracies.
For both FD and DG methods, the first order schemes produce the maximum mass flux error, spread across a larger length of the domain when compared to the third and fifth
order schemes. Between the FDLCDZW and DG methods, the spread of the mass flux error is more for the FDLCDZW schemes than that for the DG schemes but the
maximum mass flux error is slightly lower for the FDLCDZW schemes (as is evident from table \ref{tab:maxMassFluxErrors}).
\begin{table}[!htbp]
  \begin{center}
    \caption{Maximum mass flux errors across different shocks at $t=100.0$, for different schemes}
\label{tab:maxMassFluxErrors}
\begin{tabular}{|>{\centering\arraybackslash}m{0.08\textwidth}|>{\centering\arraybackslash}m{0.12\textwidth}|>{\centering\arraybackslash}m{0.15\textwidth}|>{\centering\arraybackslash}m{0.15\textwidth}|>{\centering\arraybackslash}m{0.08\textwidth}|>{\centering\arraybackslash}m{0.08\textwidth}|>{\centering\arraybackslash}m{0.08\textwidth}|}
      \hline
      \multirow{3}{5em}{Mach Number} & \multicolumn{6}{c|}{Maximum mass flux error percentage}\\\cline{2-7}
                                     &\multirow{2}{5em}{FDLR1-LF} & \multicolumn{2}{c|}{FDLCDZW-LF with FOA} & \multicolumn{3}{c|}{DG-LF}\\\cline{3-7}
                                     & & 3 & 5 & $P^0$ & $P^2$ & $P^4$ \\ \hline 
      2.0 & 14.2 & 10.8 & 9.9 & 14.1 & 12.5 & 12.5\\ \hline
      2.4 & 20.0 & 15.3 & 14.1 & 19.8 & 16.2 & 16.2 \\ \hline
      2.8 & 24.7 & 18.9 & 17.8 & 24.2 & 19.4 & 19.3 \\ \hline
      3.0 & 26.6 & 20.5 & 19.4 & 26.0 & 21.3 & 21.3 \\ \hline
    \end{tabular}
  \end{center}
\end{table}

\subsection{Quasi-One-dimensional Euler equations}\label{Sec:quasiOneD}
The quasi-one-dimensional Euler equations are given by:
\begin{equation}
  \label{eq:quasiOneDimEulerEqn}
  \frac{\partial}{\partial t}Q(x, t) + \frac{\partial}{\partial x}E(x, t) = -\frac{A^{'}(x)}{A(x)}S,
\end{equation}
where $Q = [\rho, \rho u , \rho e_t]^T$, $E(x, t) = [\rho u, \rho u^2, (\rho e_t + p)u]^T$, $S = [\rho u, \rho u^2, (\rho e_t + p)u]^T$, $A(x)$ is the area of
cross-section, with $\gamma = 1.4$,

We solve system of equations (\ref{eq:quasiOneDimEulerEqn}) in the domain $0 \leq x \leq 1$, for $A(0) = 1.0$, $A^{'}(x) = 1.0$. $x=0$ is a supersonic inflow, with
$(\rho, u, p) = (\gamma, M, 1.0)$. $x=1$ is a subsonic outflow with a back pressure $p_{back}$. The back pressure $p_{back}$ is set such that there is a shock
at $x=0.5$. Initial conditions correspond to $(\rho, u, p) = (\gamma, 0, 1)$. We obtain numerical solutions for $M=2, 2.4, 2.8, 3.0$ using FDLCDZW5-LF and DG-$P^4$-LF
methods. For 
quasi-one-dimensional Euler equations the quantity $A\rho u$ is conserved. However, in the numerical solution it is not conserved, due to using
flux splitting or a numerical flux function. For the numerical solution of one-dimensional Euler equations, using the ROE flux does not lead to any mass flux error but
it does, for quasi-one-dimensional Euler equations.

\begin{figure}[!htbp]
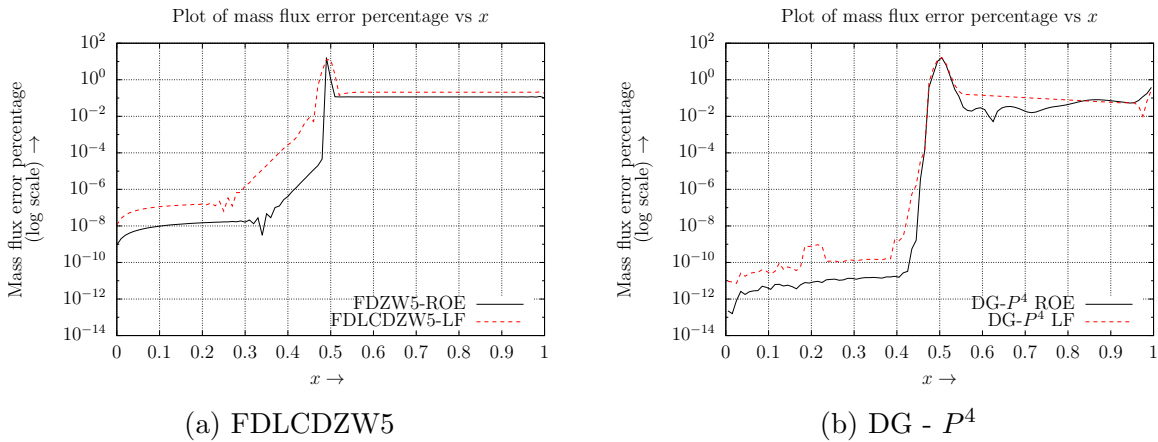

  \centering
  \begin{subfigure}{.5\textwidth}
    \centering
  \scalebox{.6}{\input{images/massFluxErr/quasi1dM2ZWenoMassFluxErr.imgtex}}
  \caption{FDLCDZW5}
  \label{fiig:quasi1dM2ZWenoMassFluxErr}
\end{subfigure}%
\begin{subfigure}{.5\textwidth}
  \centering
  \scalebox{.6}{\input{quasi1dM2DGMassFluxErr.imgtex}}
  \caption{DG - $P^4$}
  \label{fig:quasi1dM2DGP4MassFluxErr}
\end{subfigure}
\caption{PE($A(x)\rho u_i$) vs $x$, at $t=100.0$, for a quasi-one-dimensional flow with shock at $x=0.0$ with an inflow Mach Number $2.0$}
\label{fig:quasi1dM2MassFluxErr}
\end{figure}

Figure \ref{fig:quasi1dM2MassFluxErr} has plots of mass flux error for different schemes. Table \ref{tab:quasi1dmaxMassFluxErrors} has the maximum mass flux errors.
\begin{table}[!htbp]
  \begin{center}
    \caption{Maximum mass flux errors at $t=100.0$, for different inflow Mach numbers and different schemes for quasi-one-dimensional Euler equations, with shock at $x=0.5$, }
\label{tab:quasi1dmaxMassFluxErrors}
\begin{tabular}{|>{\centering\arraybackslash}m{0.1\textwidth}|>{\centering\arraybackslash}m{0.1\textwidth}|>{\centering\arraybackslash}m{0.1\textwidth}|>{\centering\arraybackslash}m{0.1\textwidth}|>{\centering\arraybackslash}m{0.1\textwidth}|}
      \hline
      \multirow{3}{5em}{Mach Number} & \multicolumn{4}{c|}{Maximum mass flux error percentage}\\\cline{2-5}
                                     &\multicolumn{2}{c|}{FDZW5} &\multicolumn{2}{c|}{DG $P^{4}$}\\\cline{2-5}
                                      & LCD LF& ROE & LF & ROE \\ \hline 
      2.0 & 12.6 & 16.1 & 16.3 & 16.1\\ \hline
      2.4 & 15.0 & 20.8 & 20.4 & 20.2 \\ \hline
      2.8 & 17.3 & 24.5 & 25.2 & 24.5 \\ \hline
      3.0 & 18.3 & 26.1 & 27.3 & 26.8 \\ \hline
    \end{tabular}
  \end{center}
\end{table}

\subsection{Two-dimensional Euler equations}
\begin{figure}[!htbp]
\begin{center}
  \includegraphics[width=0.65\textwidth]{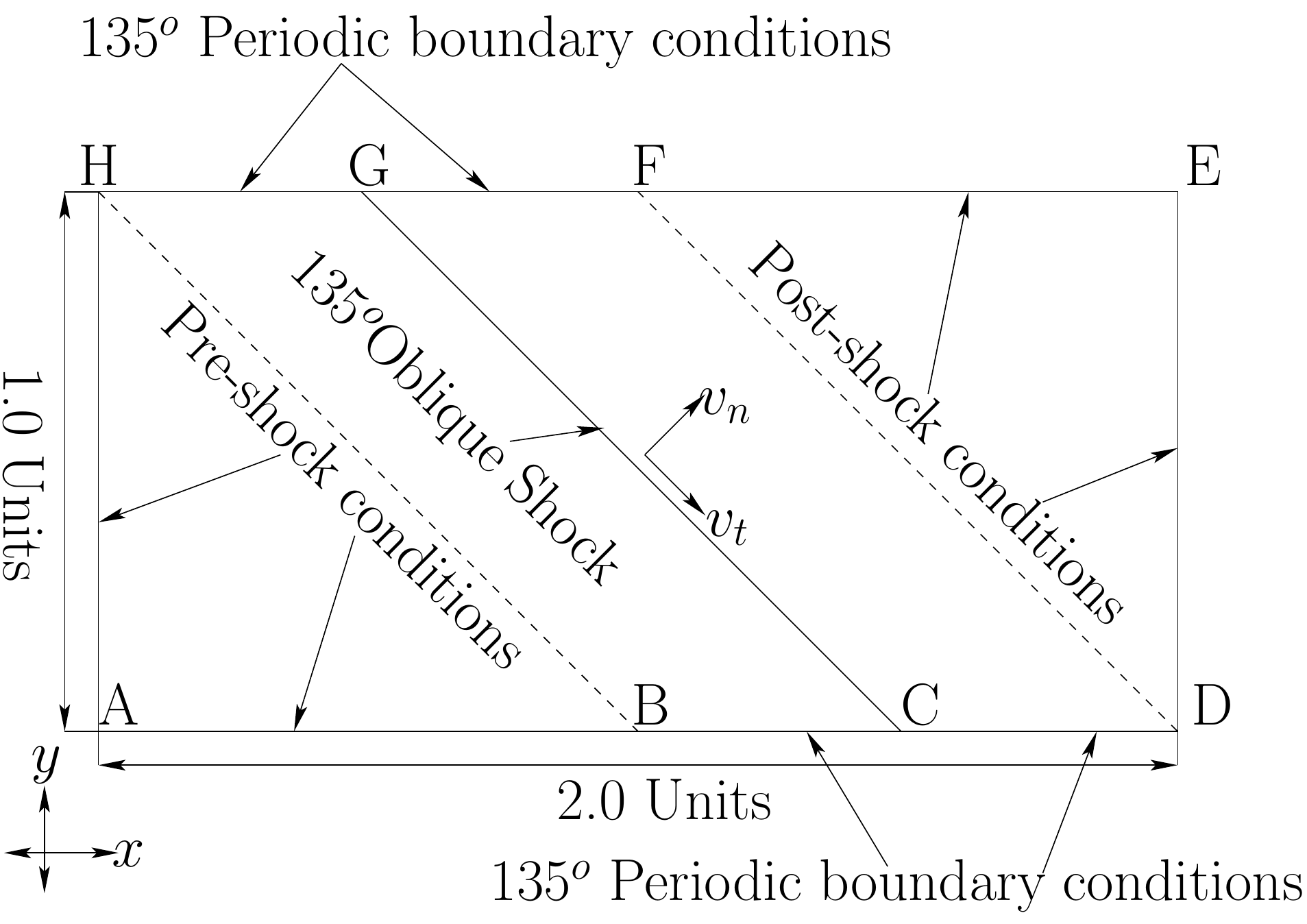}
  \caption{$135^o$ oblique shock with periodic (along $135^o$ lines) boundary conditions.}
  \label{fig:135DegOblShkGeo}
\end{center}
\end{figure}
Next, we solve the two-dimensional Euler equations (\ref{eq:Euler2dDiff}) using FDLCDZW5, FDZW5 and DG-$P^4$ with Lax-Freidrichs and ROE Fluxes. 
The computational domain
and the corresponding boundary conditions are shown in figure \ref{fig:135DegOblShkGeo}. The portion of the domain ABCGH is initialised with the `pre-shock conditions',
which are $(\rho, u, v, p) = (\gamma, M, 0, 1.0)$. The portion of the domain CDEFG is initialised with the post-shock conditions.
Let $v_n, v_t$ be components of velocity normal and parallel to the shock respectively (see figure \ref{fig:135DegOblShkGeo}). Across a stationary oblique shock 
$\rho v_n$ should be constant, but in numerical solutions obtained using Lax-Freidrichs flux and ROE flux, it is not constant.

\begin{figure}[!htbp]
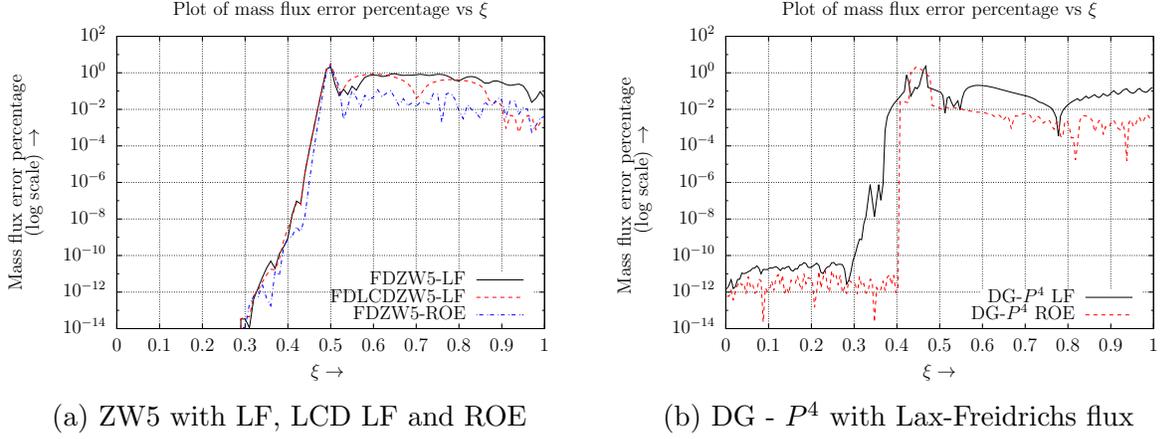

  \centering
  \begin{subfigure}{.5\textwidth}
    \centering
  \scalebox{.6}{\input{images/massFluxErr/135M2OblShkU1ZWeno5MassFluxErr.imgtex}}
  \caption{ZW5 with LF, LCD LF and ROE}
  \label{fiig:135M2OblShkU1ZWeno5MassFluxErr}
\end{subfigure}%
\begin{subfigure}{.5\textwidth}
  \centering
  \scalebox{.6}{\input{135M2OblShkDGP4MassFluxErr.imgtex}}
  \caption{DG - $P^4$ with Lax-Freidrichs flux}
  \label{fig:135M2OblShkDGP4LFMassFluxErr}
\end{subfigure}
\caption{Mass flux error percentage vs $\xi$ across a Mach 2.0 oblique shock with Lax-Freidrichs flux.}
\label{fig:135M2OblShkLF}
\end{figure}
\begin{table}[!htbp]
  \begin{center}
    \caption{Maximum mass flux errors across different $135^0$ oblique shocks for different schemes}
\label{tab:2dOblShkmaxMassFluxErrors}
\begin{tabular}{|>{\centering\arraybackslash}m{0.1\textwidth}|>{\centering\arraybackslash}m{0.1\textwidth}|>{\centering\arraybackslash}m{0.1\textwidth}|>{\centering\arraybackslash}m{0.1\textwidth}|>{\centering\arraybackslash}m{0.1\textwidth}|>{\centering\arraybackslash}m{0.1\textwidth}|}
      \hline
      \multirow{3}{5em}{Mach Number} & \multicolumn{5}{c|}{Maximum mass flux error percentage}\\\cline{2-6}
                                     & \multicolumn{3}{c|}{FDZW5} & \multicolumn{2}{c|}{DG $P^{4}$}\\\cline{2-6}
                                      & LF & LCD LF & ROE & LF & ROE \\ \hline
      2.0 & 3.5 & 3.1 & 2.9 & 2.3 & 2.03 \\ \hline
      2.4 & 7.7 & 5.0 & 5.8 & 5.2 & 5.44 \\ \hline
      2.8 & 7.0 & 8.3 & 3.4 & 6.85 & 5.92 \\ \hline
      3.0 & 6.1 & 7.6 & 5.7 & 7.4 & 6.1 \\ \hline
    \end{tabular}
  \end{center}
\end{table}
If the point A is $(0,0)$ (see figure \ref{fig:135DegOblShkGeo}), we plot the mass flux error vs a variable $\xi$ measured along the line given by $y=x-0.5$,
starting with the point $(0.5, 0)$ corresponding to $\xi=0$ and ending with $(1.5, 1)$ corresponding to $\xi=1$. 
Figure \ref{fig:135M2OblShkLF} has plots of mass flux error percentage vs $\xi$. Table \ref{tab:2dOblShkmaxMassFluxErrors} has the maximum mass flux error percentages
along the same line ($y=x-0.5$) for different schemes and different Mach numbers.

\subsection{Cause of the mass flux error}
The cause of the mass flux error seems to be the flux splitting and the approximate Riemann solver used for FDZW and DG schemes respectively. When a shock is captured
with dissipation, dissipation is also introduced in the mass conservation equation (as suggested by Jin et al \cite{jin1996}) and this leads to the mass flux error.
This mass flux error, as shown previously is high near the shock and it propagates into the region downstream of the shock.
In case the flux splitting or the approximate Riemann solver used leads to capturing the shock without any dissipation, the mass flux error is absent. This is apparent in 
the case of using Roe-flux for the normal shock problem for one-dimensional Euler equations.

\subsubsection{Alignment of shock to cell faces and Rotated Riemann flux}
For the two-dimensional Euler equations, the Roe flux also caused mass flux error because, the 135\textdegree{} shock is not aligned either with the $x$ or the $y$ 
direction. To rectify this in case of the DG method, a mesh made of right angular triangles can be chosen such that the oblique shock is aligned to the hypotenuse faces
of the triangles. Since the shock is aligned with the faces of the triangles, the Roe flux will be applied along a direction normal to the shock and therefore,
the shock is sustained without any dissipation and mass flux error.

In case of the FDZW schemes, a similar solution of applying the Roe flux in a direction normal to the shock can be used.
Let $Q_l$ and $Q_r$ be the left and right biased WENO interpolations at $(x_{i+\frac{1}{2}}, y_j)$. A rotated Riemann based flux splitting can be used with the
shock orientation being given by
\begin{multline}
  \label{eq:rotationCosineSine}
  \big(\cos(\theta), \sin(\theta)\big) = \begin{cases}
    \bigg( \frac{|u_l - u_r|}{|\vec{V}_{diff}|}, \frac{|v_l - v_r|}{|\vec{V}_{diff}|} \bigg), &\text{ if }|\vec{V}_{diff}|>10^{-2} \\
    (1, 0), &\text{ otherwise}.\\
  \end{cases},\\
  \text{where }\vec{V}_{diff} = ( |u_l - u_r|, |v_l - v_r| ),
\end{multline}
where $\theta$ is the angle made by the shock with the $y$-axis.
This will ensure that the oblique shock is captured without mass flux error and dissipation.

Capturing curved or oblique shocks or shock reflections without dissipation is a much more challenging problem as aligning the mesh faces with
shock becomes an issue for the DG method. For the FDZW method, a similar problem of shock passing through a grid point makes it challenging.
Additionally, for the FDZW method, using the rotated Riemann based flux splitting can be an issue because the procedure
given for calculating $\cos(\theta)$ and $\sin(\theta)$ in equation (\ref{eq:rotationCosineSine}) is known to lead to convergence stalling \cite{levy1993}.

Next, we look at how the introduction of viscous fluxes changes the mass flux error, with the help of the one-dimensional compressible viscous
fluid flow equations
(Newtonian fluid, Stokes' hypothesis used and with viscous and heat flux coefficients modelled using the Sutherland formulae).

\subsection{One-dimensional Compressible Viscous fluid flow equations (Navier-Stokes)}\label{Sec:compViscFluFlow}
The non-dimensional or scaled viscous fluid flow equations for Newtonian fluid in one space dimension are given by
\begin{align}
  \label{eq:viscFluidFlowEqns}
  \frac{\partial}{\partial t^*} & Q^*(x^*, t^*) + \frac{\partial}{\partial x^*} \bigg( E^*(x^*, t^*) - E^*_v(x^*, t^*) \bigg) = 0,\text{ where,}\\
  \label{eq:viscFluidFlowEqnsTerms1}
Q^* = &\begin{bmatrix} \rho^* \\ \rho^* u^* \\ \rho^* e^*_t \end{bmatrix},
 E^* = \begin{bmatrix} \rho^* u^* \\ \rho^* u^{*2} \\ (\rho^* e^*_t + p^*)u^* \end{bmatrix},
    e^{*}_t = \frac{p^{*}}{\rho^{*}(\gamma -1)} + \frac{1}{2}\left(u^{*2} + v^{*2} \right),\\
  \label{eq:viscFluidFlowEqnsTerms2}
    E^*_v = &\bigg[ 0~~~~ \left(\frac{\partial u^*}{\partial x^*}\frac{\lambda+2\mu}{\rho_0 U_0 L}\right)~~~~
  \left(u^*\frac{\partial u^*}{\partial x^*}\frac{\lambda+2\mu}{\rho_0 U_0 L} + \frac{\kappa}{R\rho_0U_0L}\frac{\partial T^*}{\partial x^*}\right) \bigg]^T
\end{align}
The scaling used is:
\begin{multline}
  \label{eq:viscFluidFlowEqnsTerms3}
p = p^*\rho_0~ U_0^2,~ e_t = e_t^* U_0^2, x = x^* L,~ u = u^* U_0,~ t = t^*\frac{L}{U_0},~ \rho = \rho^*~ \rho_0,
  a = a^* U_0,\\\text{ and } T = T^*\frac{U_0^2}{R},
\end{multline}
where, $R = 287.4J/(kg K), \gamma = 1.4$. The Stokes' hypothesis is assumed which is $3\lambda + 2\mu = 0$. The Sutherland model for coefficients of viscosity and
heat conduction is used,
given by
\begin{multline}
  \label{eq:sutherlandViscModel}
  \mu = C_1\frac{T^{\frac{3}{2}}}{T+C_2},~\kappa = C_3\frac{T^{\frac{3}{2}}}{T+C_4},\text{ where }C_1 = 1.458X10^{-6} \frac{kg}{ms\sqrt{K}}, C_2 = 110.4K,\\
  C_3 = 2.495X10^{-3}\frac{kgm}{s^3K^{\frac{3}{2}}},\text{ and }C_4 = 194K.
\end{multline}
The values of the scaling parameters are $\rho_0 = 1.204kg/m^3, U_0 = 343.249m/s$.
System of equations (\ref{eq:viscFluidFlowEqns}) are solved in the domain $0 \leq x^{*} \leq 1$ with initial conditions
\begin{equation}
  \label{eq:oneDNormShkVisc}
  Q^*(x^*, 0) = \begin{cases}
    Q^*_{BS} & x^*<0.5 \\
    Q^*_{AS} & x^*\geq 0.5\\
  \end{cases},
\end{equation}
where
\begin{equation}
  \label{eq:1dNormShkViscInitConds}
\begin{bmatrix}\rho^*_{BS} \\ u^*_{BS} \\ p^*_{BS}\end{bmatrix} = \begin{bmatrix}\gamma \\ M \\ 1.0\end{bmatrix},
\begin{bmatrix}\rho^*_{AS} \\ u^*_{AS} \\ p^*_{AS}\end{bmatrix} = 
\begin{bmatrix}\frac{(\gamma+1)M^2 \rho^*_{BS}}{(\gamma-1)M^2 + 2}
\\ \frac{\rho^*_{BS}u^*_{BS}}{\rho^*_{AS}} \\ 
\frac{p^*_{BS}(2\gamma M^2 - (\gamma -1))}{\gamma+1}
\end{bmatrix},
\end{equation}
with supersonic inflow conditions at $x^* = 0.0$ and subsonic outflow conditions with back pressure of $p^*_{AS}$ at $x^* = 1.0$. We chose a mesh with GPS( $ = \Delta x^*$)
of $1/100$. We use the numerical methods described in section \ref{Sec:NumMethod} for obtaining the numerical solutions. As mentioned in section \ref{Sec:extnToViscEqns}
we use FDLR5 (equation \ref{eq:linHiOrd}, $\hat{A} = 0$) to calculate viscous fluxes and their derivatives for the FD schemes. For calculating inviscid fluxes and their
derivatives,
we use FDZW5-LF or FDLCDZW5-LF or FDZW5-ROE. Also, we obtain numerical solutions using LDG method.
We obtain numerical solutions for inflow conditions corresponding to $M = 2.0, 2.4, 2.8, 3.0$ for the values of parameter $L$ equal to $1.0, 10^{-4},\text{ and } 10^{-6}$.

For the viscous fluid flow equations, using Roe flux also leads to mass flux error. Table \ref{tab:maxMassFluxErrorsVisc} has the maximum mass flux error percentages for
different schemes for $L=1$ and $10^{-4}$. The maximum mass flux error for $L=10^{-6}$ using FDZW5-LF, FDLCDZW5-LF, FDZW5-ROE, and LDG-P4 with Lax-Freidrichs or ROE
Flux are of the order of $10^{-3}$.
\begin{table}[!htbp]
  \begin{center}
    \caption{Maximum mass flux errors across different shocks at $t=100.0$, for different schemes}
\label{tab:maxMassFluxErrorsVisc}
\begin{tabular}{|>{\centering\arraybackslash}m{0.1\textwidth}|>{\centering\arraybackslash}m{0.1\textwidth}|>{\centering\arraybackslash}m{0.1\textwidth}|>{\centering\arraybackslash}m{0.1\textwidth}|>{\centering\arraybackslash}m{0.1\textwidth}|>{\centering\arraybackslash}m{0.1\textwidth}|}
      \hline
      \multirow{3}{5em}{Mach Number} & \multicolumn{5}{c|}{Maximum mass flux error percentage for $L=1$}\\\cline{2-6}
                                     &\multicolumn{3}{c|}{FDZW5} &\multicolumn{2}{c|}{LDG-$P^4$}\\\cline{2-6}
                                     & LF & LCD LF & ROE & LF & ROE\\ \hline 
      2.0 & 13.4 & 9.7 & 2.2 & 9.9 & 2.03\\ \hline
      2.4 & 18.6 & 13.8 & 9.1 & 14.1 & 8.7 \\ \hline
      2.8 & 22.7 & 17.1 & 16.0 & 17.5 & 15.8 \\ \hline
      3.0 & 24.4 & 18.5 & 19.1 & 18.2 & 19.3 \\ \hline
      \multicolumn{6}{|c|}{Maximum mass flux error percentage for $L=10^{-4}$}\\\cline{1-6}
      2.0 & 13.2 & 9.2 & 10.4 & 9.7 & 10.3\\ \hline
      2.4 & 18.6 & 13.12 & 14.4 & 13.4 & 13.9 \\ \hline
      2.8 & 21.9 & 16.4 & 23.2 & 16.7 & 16.2 \\ \hline
      3.0 & 23.6 & 17.8 & 24.4 & 18.3 & 18.1 \\ \hline
    \end{tabular}
  \end{center}
\end{table}
\begin{figure}[!htbp]
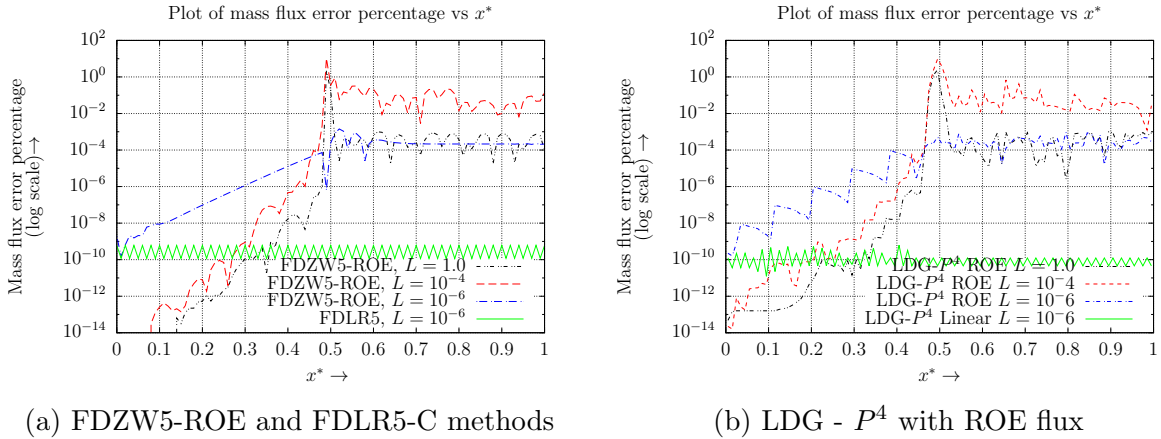

  \centering
  \begin{subfigure}{.5\textwidth}
    \centering
  \scalebox{.6}{\input{images/massFluxErr/Mach2ViscROEMassFluxErr.imgtex}}
  \caption{FDZW5-ROE and FDLR5-C methods}
  \label{fig:massFluxErrM2ZW5ROEVisc}
\end{subfigure}%
\begin{subfigure}{.5\textwidth}
  \centering
  \scalebox{.6}{\input{Mach2ViscROEMassFluxErrDG.imgtex}}
  \caption{LDG - $P^4$ with ROE flux}
  \label{fig:massFluxErrM2DGP4ROEVisc}
\end{subfigure}
\caption{Mass flux error percentage vs $x^*$ at $t=100.0$ across a Mach 2.0 shock using ROE flux for WENO and DG schemes for different values of $L$.}
\label{fig:massFluxErrM2Visc}
\end{figure}

For $L=10^{-6}$ and $\Delta x^* = 1/100$ the shock is sufficiently resolved and therefore it is not necessary to use WENO reconstruction or limiter for the LDG method.
Also, the upwind biasing of Inviscid fluxes is not necessary and therefore one can calculate inviscid flux derivatives with $\hat{A}=0$ (see equation
\ref{eq:FSplitting}). Therefore, we use FDLR5-C and DG $P^4$-C methods . Figure \ref{fig:massFluxErrM2Visc} has plots of mass
flux error percentage vs $x^*$ for different schemes.
Clearly, the mass flux error is minimum for the numerical solution obtained using FDLR5-C or the DG $P^{4}$-C Linear schemes. 
\section{Refinement near the shock using overset meshes}\label{Sec:refinementAndOverset}
Of course using a fine mesh corresponding to a $\Delta x^{*} = 10^{-8}$ for problems of general interest may not be possible. One solution in such cases is to reduce the
mass conservation error by doing a mesh
refinement near the shock using a series of overset meshes. This is demonstrated by applying it in computing numerical solutions of quasi-one-dimensional Euler
equations and two-dimensional Euler equations.
\subsection{Quasi-One-dimensional Euler Equations}\label{Sec:quasiOneDOverset}
We solve the same problem mentioned in section \ref{Sec:quasiOneD} for $M=3.0$, using FDLCDZW5-LF, DG $P^4$-LF with TVD-RK3 method using three mesh configurations. The first
configuration (Config1) is a mesh in the domain $0 \leq x \leq 1$ with GPS or $\Delta x$ of $1/200$. The second configuration (Config2) is that with with a mesh in
the domain $0.48 \leq x \leq 0.525$, with GPS or $\Delta x$ of $1/2200$, overset on mesh config1. The third configuration (Config3) is that with a mesh in the domain
$0.499\overline{54} \leq x \leq 0.50\overline{36}$, with GPS or $\Delta x$ of $1/24200$, overset on mesh config2. First a solution is obtained using Config1. Initial
conditions for Config2 is the
numerical solution obtained using Config1 and that for Config3 is numerical solution obtained using Config2 ( See section \ref{Sec:OversetDesc} for more details ).

\begin{figure}[!htbp]
  \centering
  \begin{subfigure}{.5\textwidth}
    \centering
  \scalebox{.6}{\input{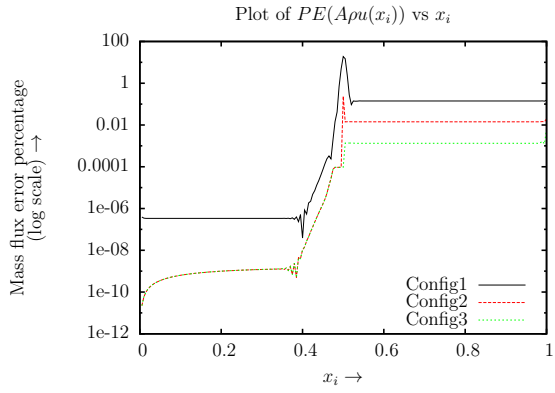}}
  \caption{PE($A\rho u(x_i)$) vs $x_i$}
  \label{fig:PEADUVsx3Configs}
\end{subfigure}%
\begin{subfigure}{.5\textwidth}
  \centering
  \scalebox{.6}{\input{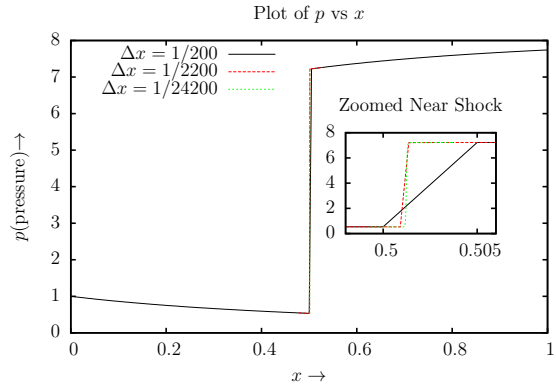}}
  \caption{pressure vs $x$ in the three meshes obtained using Config3 }
  \label{fig:PVsX3Configs}
\end{subfigure}
\caption{Plots of PE($A\rho u(x_i)$) vs $x_i$ and pressure vs $x$ for meshes Config1, Config2, Config3, obtained using FDLCDZW5-LF. }
\label{fig:oversetWENOPlots}
\end{figure}

\begin{figure}[!htbp]
  \centering
  \begin{subfigure}{.5\textwidth}
    \centering
  \scalebox{.6}{\input{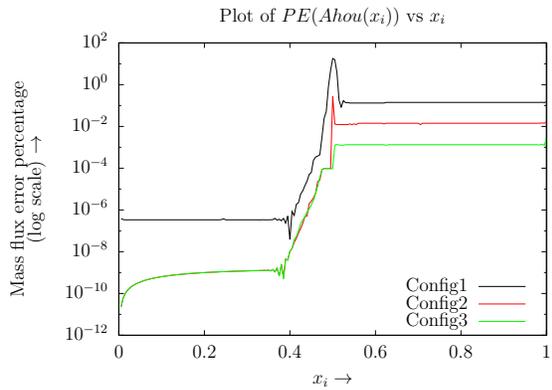}}
  \caption{PE($A\rho u(x_i)$) vs $x_i$}
  \label{fig:PEADUVsx3ConfigsDG}
\end{subfigure}%
\begin{subfigure}{.5\textwidth}
  \centering
  \scalebox{.6}{\input{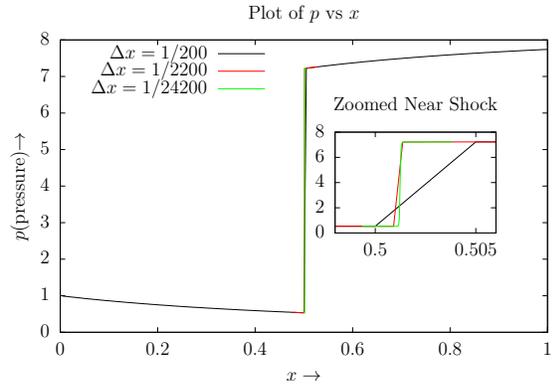}}
  \caption{pressure vs $x$ in the three meshes obtained using Config3 }
  \label{fig:PVsX3ConfigsDG}
\end{subfigure}
\caption{Plots of PE($A\rho u(x_i)$) vs $x_i$ and pressure vs $x$ for meshes Config1, Config2, Config3, obtained using DG $P^4$-LF. }
\label{fig:oversetDGPlots}
\end{figure}

As mentioned in section \ref{Sec:quasiOneD}, $\rho A u$ should be constant along the domain but it is not in the numerical solution. Figures \ref{fig:PEADUVsx3Configs},
\ref{fig:PEADUVsx3ConfigsDG} have  plots of PE($A\rho u(x_i)$) vs $x_i$ for the three mesh configurations obtained using FDLCDZW5-LF and DG $P^4$-LF methods.
Figures \ref{fig:PVsX3Configs}, \ref{fig:PVsX3ConfigsDG} have corresponding plots of pressure vs x for the three mesh configurations. 

As can be seen in the figures \ref{fig:oversetWENOPlots} and \ref{fig:oversetDGPlots}, going from Config1 to Config3,
the difference between the post shock mass flux error percentage drops from approximately $10^{-1}$ to $10^{-3}$. Therefore, the error that is propagated into the region downstream of the shock reduces by doing a mesh refinement near the shock. Also, the solutions
obtained using WENO and DG methods are almost the same.

\subsection{Curing the carbuncle: Initial results}

\begin{figure}[!htbp]
  \centering
  \begin{subfigure}{.5\textwidth}
    \centering
  \includegraphics[width=0.9\textwidth]{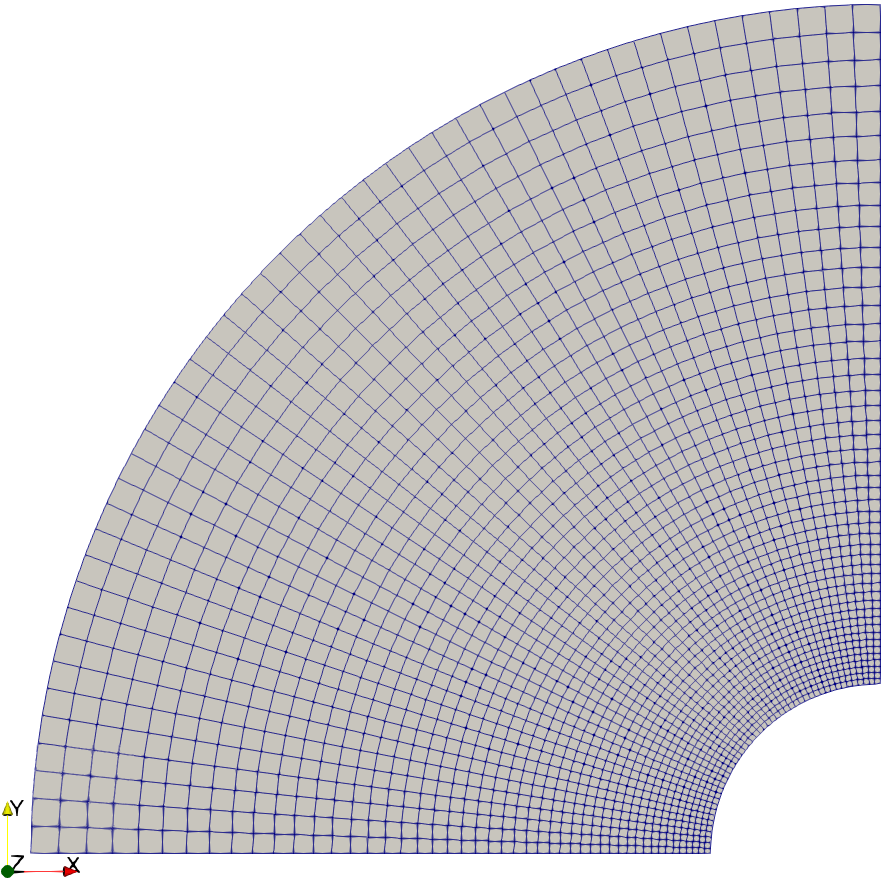}
  \caption{Structured mesh}
  \label{fig:structuredMeshFlowOverCyl}
\end{subfigure}%
\begin{subfigure}{.5\textwidth}
  \centering
  \includegraphics[width=0.9\textwidth]{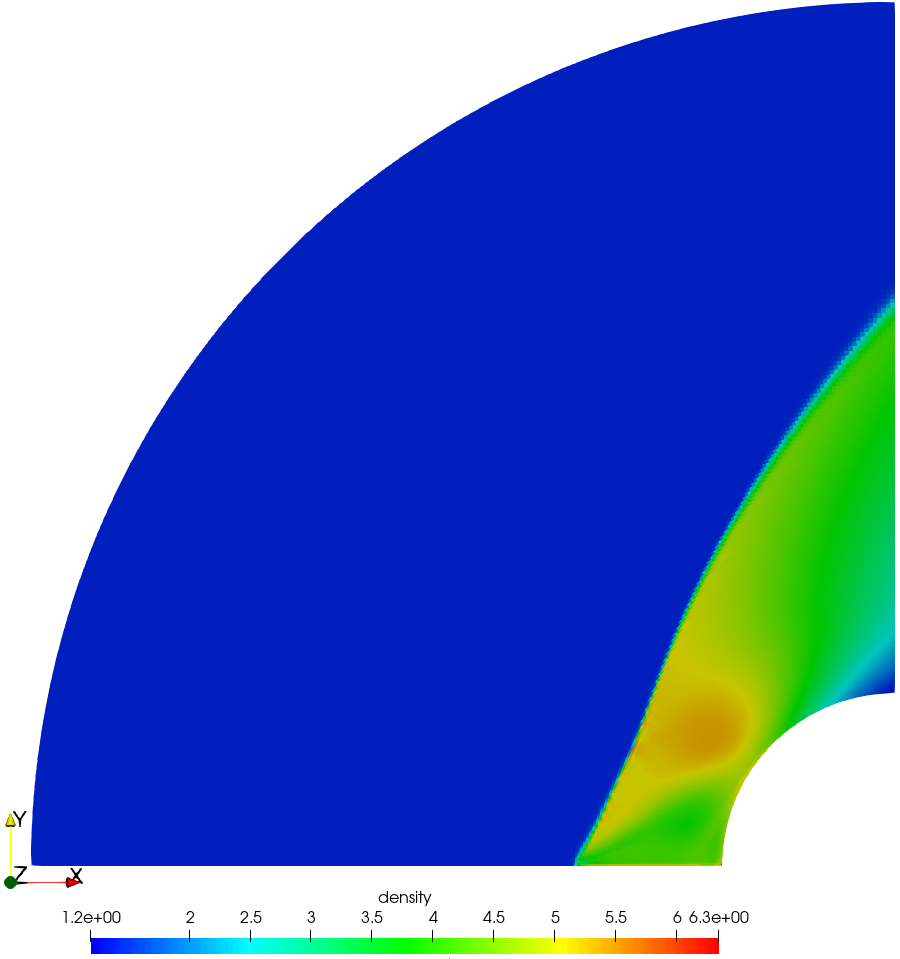}
  \caption{Colour plot of density showing Carbuncle}
  \label{fig:densityPlotWithCarbuncle}
\end{subfigure}
\caption{Mach 3.0 flow over a circular cylinder}
\label{fig:mach3FlowOverCylStrMshCarb}
\end{figure}
\begin{figure}[!htbp]
  \centering
  \begin{subfigure}{.5\textwidth}
    \centering
  \includegraphics[width=0.9\textwidth]{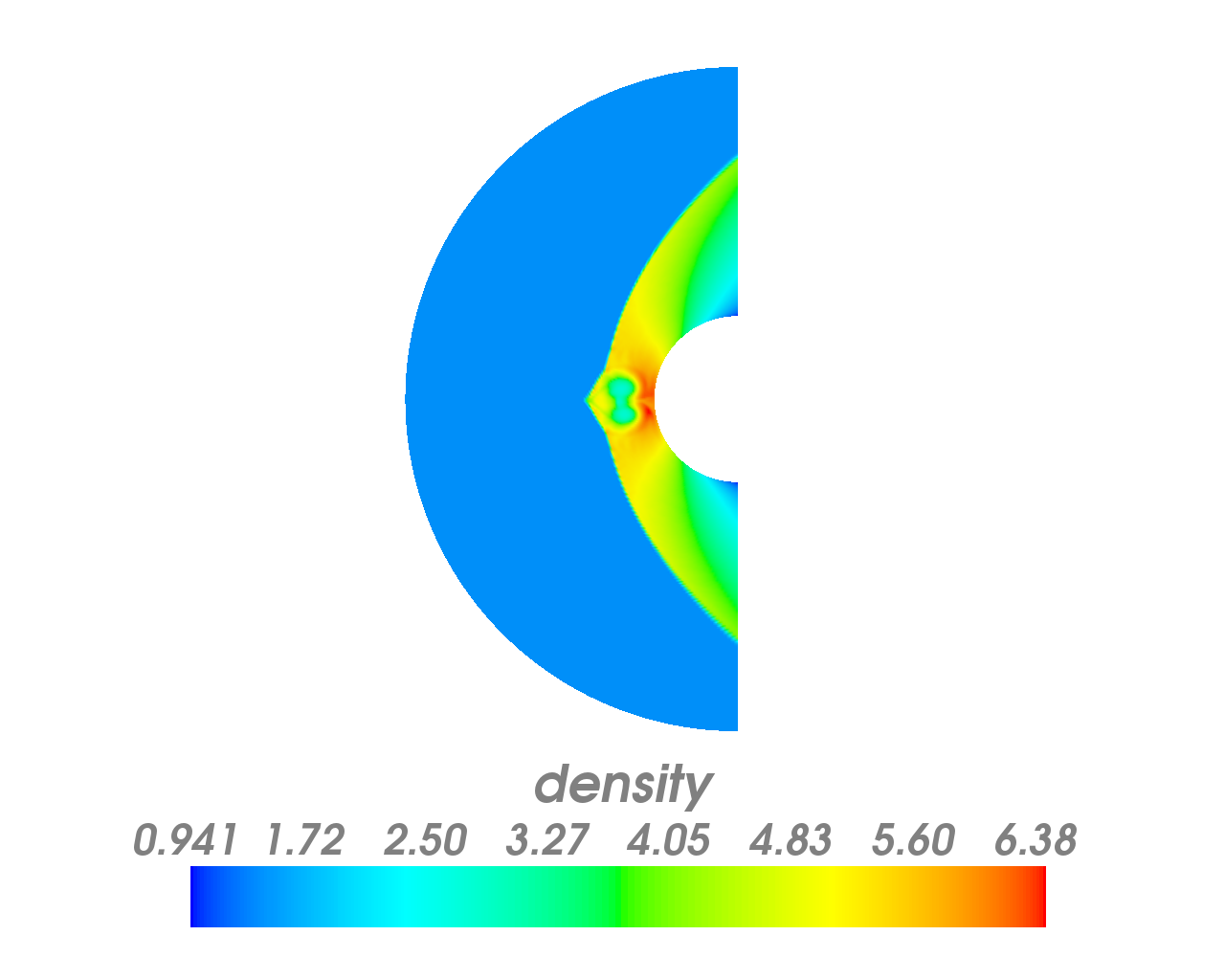}
  \caption{Density Solution with Carbuncle, without refinement}
  \label{fig:RoeFluxWCarbuncle}
\end{subfigure}%
\begin{subfigure}{.5\textwidth}
  \centering
  \includegraphics[width=0.9\textwidth]{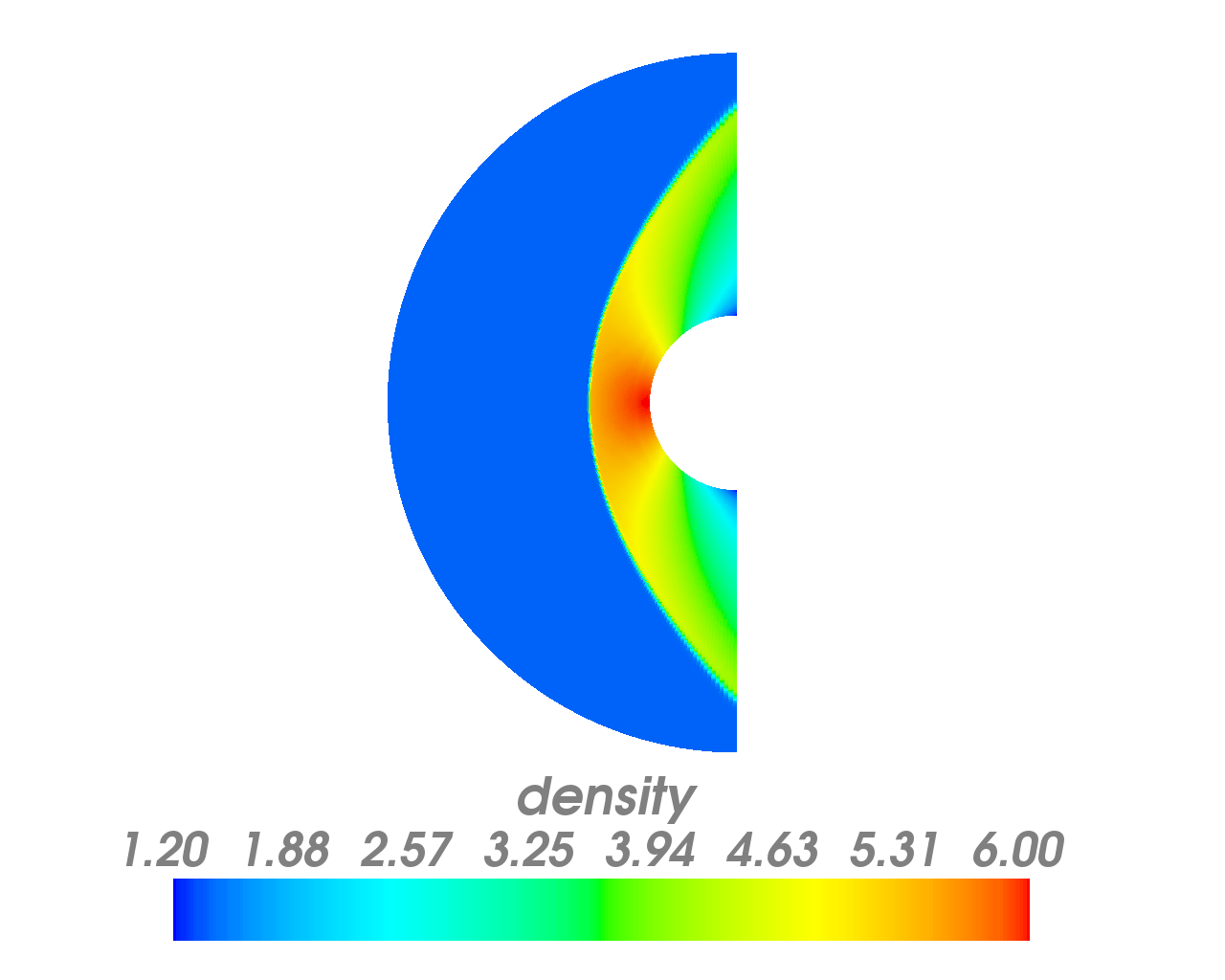}
  \caption{Density solution without carbuncle obtained using two levels of overset mesh}
  \label{fig:RoeFluxWoCarbuncle}
\end{subfigure}
\caption{Density plots for flow over a circular cylinder using Roe Flux, obtained using DG $P^2$}
\label{fig:RoeFlux}
\end{figure}
It is well known that a carbuncle (as shown in Figure \ref{fig:densityPlotWithCarbuncle}) forms when the Roe flux and a structured mesh (as shown in Figure \ref{fig:structuredMeshFlowOverCyl}) are used for computing flow over a circular cylinder. 
An often proposed solution to this problem is using dissipative
flux functions in the direction normal to the shock \cite{nishikawa2008}. The reasons for the formation of carbuncle were studied and that there may be a connection between numerical mass flux and the carbuncle formation was reported in literature \cite{liou2000}.

One way of reducing error in the
numerical mass flux, as demonstrated above, is the use of multiple overset meshes. This kind of mesh refinement near the shock, using multiple overset meshes (whilst using
Roe flux) seems to cure the carbuncle problem. We report here,
preliminary results obtained using such multiple overset mesh configuration (two levels for the present case, similar to the one described in section
\ref{Sec:quasiOneDOverset} ) in the shock region for the problem of Mach 3.0
flow over a circular cylinder. We have observed that the carbuncle disappears by using two levels of overset mesh. The solutions with and without the carbuncle,
obtained using meshes without and with two levels of overset meshes respectively, 
are shown in Figure \ref{fig:RoeFlux}. More details of this methodology will be given in a subsequent paper, focused
on accurate shock capturing using high-order methods.

\begin{figure}[!htbp]
  \centering
  \begin{subfigure}{.5\textwidth}
    \centering
  \includegraphics[width=0.9\textwidth]{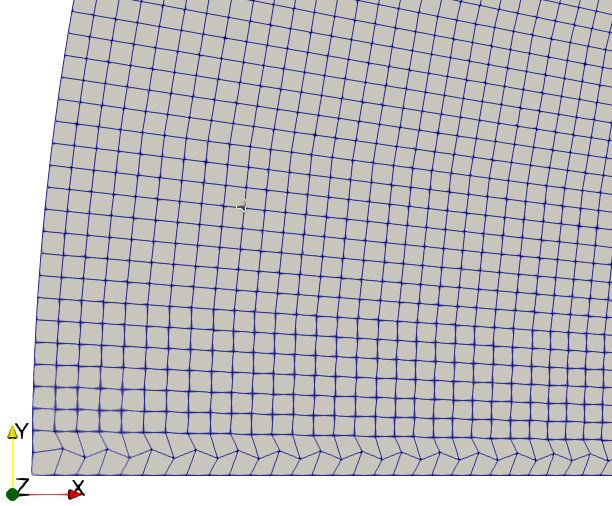}
  \caption{Carbuncle Avoiding Structured mesh: Perturbing bottom row of cells}
  \label{fig:perturbedMesh}
\end{subfigure}%
\begin{subfigure}{.5\textwidth}
  \centering
  \includegraphics[width=0.9\textwidth]{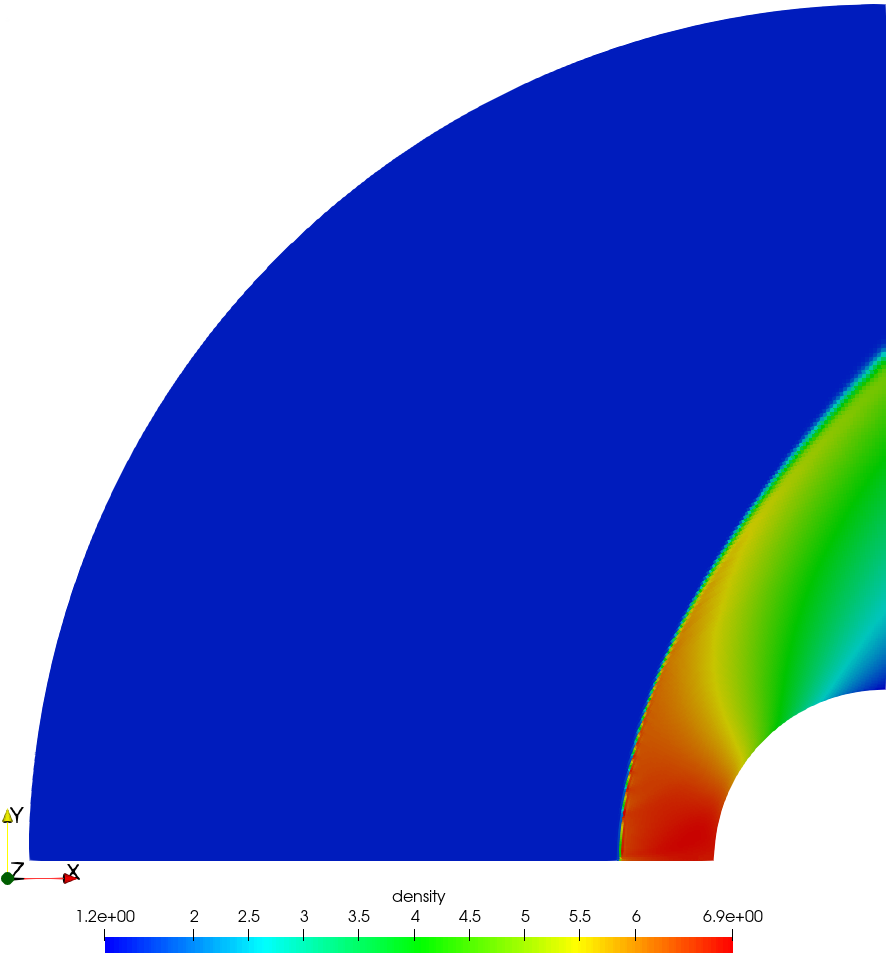}
  \caption{Density solution without carbuncle obtained using mesh similar to the one on the left}
  \label{fig:perturbedMeshMach3Soln}
\end{subfigure}
\caption{Avoiding the carbuncle by perturbing the bottom row of cells}
\label{fig:perturbedMeshAndSoln}
\end{figure}
Interestingly, another way to cure or avoid the carbuncle is to perturb the mesh such that the faces of the two bottom most row of cells are not parallel to the x and y directions as show in figure \ref{fig:perturbedMesh}. Using this mesh with the Roe flux does not produce
the carbuncle as shown in figure \ref{fig:perturbedMeshMach3Soln}.
In this mesh, the bottom portion of the shock, or the ``normal portion''( portion of the shock that is almost a normal shock) of the shock , that forms and travels upstream is not aligned with the cell faces and hence using the Roe flux also leads to introducing dissipation or mass
conservation error. Where as for a mesh similar to the one shown in figure \ref{fig:structuredMeshFlowOverCyl}, the ``normal portion'' of the shock is almost aligned with the cell faces and hence there will be essentially no mass conservation error in this
region,
but near the ``oblique portion'' of the shock (region excluding the ``normal portion''), there will be considerable mass conservation error. These results show the link between mass conservation error and carbuncle.
It shows that essentially zero mass conservation error near ``normal portion'' of the shock and considerable error near ``oblique portion'' of the shock could be the reason for
formation of Carbuncle. Also, the carbuncle can be cured or avoided
either by reducing the mass conservation error near the ``oblique portion'' of the shock by using overset mesh with 2 levels of refinement or by introducing mass conservation error near the ``normal portion'' of the shock by skewing the mesh.
\section{Conclusion}\label{Sec:Conclusions}
For the moving shock problem, we compared the performance of different numerical fluxes used in combination with different numerical methods. For the first order methods,
we showed that for certain problem parameters, ROE flux performs better than Osher flux, as opposed to the problems generally reported and cited in literature
\cite{roberts1990}. We underscored the importance of doing a characteristic-wise reconstruction for high-order methods
by giving an example of a case where doing a component-wise reconstruction instead, leads to `NAN's in the computation.

Using the test problems of normal shock for one-dimensional and quasi-one-dimensional Euler equations and the test problem of oblique shock for two-dimensional Euler equations,
we have shown that mass flux error occurs due to the use of dissipative flux splittings for conservative finite difference WENO schemes and the use of dissipative
flux functions (approximate Riemann solvers) for the Discontinuous Galerkin method. We showed that the mass flux error varies with Mach number
before the shock and formal order of accuracy of the scheme. We showed that using ROE flux also leads to significant mass conservation error while solving
the quasi-one-dimensional Euler equations, one-dimensional viscous fluid flow equations and
two-dimensional Euler equations. 

For the two-dimensional Euler equations, for the simple problem of the 135\textdegree{} oblique shock, techniques like, choosing a
mesh with shock aligned to cell faces for DG method and using the flux splitting based on Rotated Riemann solvers for the FDLCDZW, to avoid
the mass flux error were given. However, extending these techniques for capturing more complex flows having curved shocks or shock reflections is not straightforward.

We showed that without upwind biasing, using high order linear reconstruction for the conservative finite difference scheme and using a central flux function for
the DG method, a shock can be captured, if a mesh of sufficient resolution is used.

We applied the technique of using multiple levels of overset meshes for resolving flow near the shock for a quasi-one-dimensional flow problem.
We showed that using such a mesh leads to mitigation of mass flux error.
We also showed that the connection between the mass conservation error and the formation of carbuncle. We showed preliminary results of two ways of curing the carbuncle. One, by reducing the mass conservation near the ``oblique portion'' of the shock by using overset mesh with 2 
levels of refinement and the other by
introducing mass conservation error near the ``normal portion'' of the shock by skewing the mesh near the ``normal portion'' of the shock.

Finally, the mass conservation error or the post shock oscillation error seemed essentially independent of whether the WENO method was used or DG method was used.
\bibliography{research,references}

\begin{thebibliography}{45}
\expandafter\ifx\csname natexlab\endcsname\relax\def\natexlab#1{#1}\fi
\providecommand{\url}[1]{\texttt{#1}}
\providecommand{\href}[2]{#2}
\providecommand{\path}[1]{#1}
\providecommand{\DOIprefix}{doi:}
\providecommand{\ArXivprefix}{arXiv:}
\providecommand{\URLprefix}{URL: }
\providecommand{\Pubmedprefix}{pmid:}
\providecommand{\doi}[1]{\href{http://dx.doi.org/#1}{\path{#1}}}
\providecommand{\Pubmed}[1]{\href{pmid:#1}{\path{#1}}}
\providecommand{\bibinfo}[2]{#2}
\ifx\xfnm\relax \def\xfnm[#1]{\unskip,\space#1}\fi
\bibitem[{Arora and Roe(1997)}]{arora1997}
\bibinfo{author}{Arora, M.}, \bibinfo{author}{Roe, P.L.}, \bibinfo{year}{1997}.
\newblock \bibinfo{title}{On postshock oscillations due to shock capturing
  schemes in unsteady flows}.
\newblock \bibinfo{journal}{Journal of Computational Physics}
  \bibinfo{volume}{130}, \bibinfo{pages}{25 -- 40}.
\newblock \URLprefix
  \url{http://www.sciencedirect.com/science/article/pii/S0021999196955345},
  \DOIprefix\doi{https://doi.org/10.1006/jcph.1996.5534}.
\bibitem[{Barth(1989)}]{barth1989}
\bibinfo{author}{Barth, T.J.}, \bibinfo{year}{1989}.
\newblock \bibinfo{title}{Some notes on shock resolving flux functions. part 1:
  Stationary characteristics}.
\newblock \bibinfo{type}{Technical Report} \bibinfo{number}{NASA-TM-101087,
  A-89087, NAS 1.15:101087}. NASA.
\bibitem[{Borges et~al.(2008)Borges, Carmona, Costa and Don}]{borges2008}
\bibinfo{author}{Borges, R.}, \bibinfo{author}{Carmona, M.},
  \bibinfo{author}{Costa, B.}, \bibinfo{author}{Don, W.S.},
  \bibinfo{year}{2008}.
\newblock \bibinfo{title}{An improved weighted essentially non-oscillatory
  scheme for hyperbolic conservation laws}.
\newblock \bibinfo{journal}{Journal of Computational Physics}
  \bibinfo{volume}{227}, \bibinfo{pages}{3191 -- 3211}.
\newblock \URLprefix
  \url{http://www.sciencedirect.com/science/article/pii/S0021999107005232},
  \DOIprefix\doi{https://doi.org/10.1016/j.jcp.2007.11.038}.
\bibitem[{Castro et~al.(2011)Castro, Costa and Don}]{castro2011}
\bibinfo{author}{Castro, M.}, \bibinfo{author}{Costa, B.},
  \bibinfo{author}{Don, W.S.}, \bibinfo{year}{2011}.
\newblock \bibinfo{title}{High order weighted essentially non-oscillatory
  weno-z schemes for hyperbolic conservation laws}.
\newblock \bibinfo{journal}{Journal of Computational Physics}
  \bibinfo{volume}{230}, \bibinfo{pages}{1766 -- 1792}.
\newblock \URLprefix
  \url{http://www.sciencedirect.com/science/article/pii/S0021999110006431},
  \DOIprefix\doi{https://doi.org/10.1016/j.jcp.2010.11.028}.
\bibitem[{Cheng et~al.(2013)Cheng, Lu and Liu}]{cheng2013}
\bibinfo{author}{Cheng, J.}, \bibinfo{author}{Lu, Y.}, \bibinfo{author}{Liu,
  T.}, \bibinfo{year}{2013}.
\newblock \bibinfo{title}{Multidomain hybrid rkdg and weno methods for
  hyperbolic conservation laws}.
\newblock \bibinfo{journal}{SIAM Journal on Scientific Computing}
  \bibinfo{volume}{35}, \bibinfo{pages}{A1049--A1072}.
\newblock \URLprefix \url{https://doi.org/10.1137/110855156},
  \DOIprefix\doi{10.1137/110855156},
  \href{http://arxiv.org/abs/https://doi.org/10.1137/110855156}{{\tt
  arXiv:https://doi.org/10.1137/110855156}}.
\bibitem[{Cheng et~al.(2016)Cheng, Wang and Liu}]{cheng2016}
\bibinfo{author}{Cheng, J.}, \bibinfo{author}{Wang, K.}, \bibinfo{author}{Liu,
  T.}, \bibinfo{year}{2016}.
\newblock \bibinfo{title}{A general high-order multi-domain hybrid dg/weno-fd
  method for hyperbolic conservation laws}.
\newblock \bibinfo{journal}{Journal of Computational Mathematics}
  \bibinfo{volume}{34}, \bibinfo{pages}{30--48}.
\bibitem[{Cockburn et~al.(1990)Cockburn, Hou and Shu}]{chs}
\bibinfo{author}{Cockburn, B.}, \bibinfo{author}{Hou, S.},
  \bibinfo{author}{Shu, C.}, \bibinfo{year}{1990}.
\newblock \bibinfo{title}{{TVB} {Runge-Kutta} local projection discontinuous
  {Galerkin} finite element method for scalar conservation laws {IV:} {The}
  {multi-dimensional} case.}
\newblock \bibinfo{journal}{Math. Comp.} \bibinfo{volume}{54},
  \bibinfo{pages}{545--581}.
\bibitem[{Cockburn et~al.(2000)Cockburn, Karniadakis and Shu}]{cks}
\bibinfo{author}{Cockburn, B.}, \bibinfo{author}{Karniadakis, G.},
  \bibinfo{author}{Shu, C.W.}, \bibinfo{year}{2000}.
\newblock \bibinfo{title}{{Discontinuous} {Galerkin} {Methods:} {Theory,}
  {Computation} and {Applications}}. \bibinfo{publisher}{Lecture {Notes} in
  {Computational} {Science} and {Engineering,} {Springer,} {vol:11}}. chapter
  \bibinfo{chapter}{The development of discontinuous {Galerkin} {Methods}}.
\newblock pp. \bibinfo{pages}{3--50}.
\bibitem[{Cockburn and Shu(1989)}]{cs2}
\bibinfo{author}{Cockburn, B.}, \bibinfo{author}{Shu, C.},
  \bibinfo{year}{1989}.
\newblock \bibinfo{title}{{TVB} {Runge-Kutta} local projection discontinuous
  {Galerkin} finite element method for scalar conservation laws {II}: General
  frame work}.
\newblock \bibinfo{journal}{Math. Comp.} \bibinfo{volume}{52},
  \bibinfo{pages}{411--435}.
\bibitem[{Cockburn and Shu(1998a)}]{cs3}
\bibinfo{author}{Cockburn, B.}, \bibinfo{author}{Shu, C.},
  \bibinfo{year}{1998}a.
\newblock \bibinfo{title}{The {Runge-Kutta} discontinuous {Galerkin} finite
  element method for conservation laws {V:} {Multidimensional} systems.}
\newblock \bibinfo{journal}{Journal of Computational Physics}
  \bibinfo{volume}{141}, \bibinfo{pages}{199--224}.
\bibitem[{Cockburn and Shu(1991)}]{cs1}
\bibinfo{author}{Cockburn, B.}, \bibinfo{author}{Shu, C.W.},
  \bibinfo{year}{1991}.
\newblock \bibinfo{title}{The runge kutta local projection
  {$P^{1}$-discontinuous} {Galerkin} method for scalar conservation laws.}
\newblock \bibinfo{journal}{{$M^{2}AN$}} \bibinfo{volume}{25},
  \bibinfo{pages}{337--361}.
\bibitem[{Cockburn and Shu(1998b)}]{cs5}
\bibinfo{author}{Cockburn, B.}, \bibinfo{author}{Shu, C.W.},
  \bibinfo{year}{1998}b.
\newblock \bibinfo{title}{The local {Discontinuous Galerkin Method} for
  {time-dependent convection-diffusion} systems.}
\newblock \bibinfo{journal}{SIAM J. Numer. Anal.} \bibinfo{volume}{35},
  \bibinfo{pages}{2440--2463}.
\bibitem[{Galbraith et~al.(2014)Galbraith, Benek, Orkwis and Turner}]{gbot}
\bibinfo{author}{Galbraith, M.}, \bibinfo{author}{Benek, J.},
  \bibinfo{author}{Orkwis, P.}, \bibinfo{author}{Turner, M.},
  \bibinfo{year}{2014}.
\newblock \bibinfo{title}{A {Discontinuous Galerkin Chimera scheme}.}
\newblock \bibinfo{journal}{Computers and Fluids} \bibinfo{volume}{98},
  \bibinfo{pages}{27--53}.
\bibitem[{Gottlieb et~al.(2001)Gottlieb, Shu and Tadmor}]{gottlieb2001}
\bibinfo{author}{Gottlieb, S.}, \bibinfo{author}{Shu, C.W.},
  \bibinfo{author}{Tadmor, E.}, \bibinfo{year}{2001}.
\newblock \bibinfo{title}{Strong stability-preserving high-order time
  discretization methods}.
\newblock \bibinfo{journal}{SIAM Review} \bibinfo{volume}{43},
  \bibinfo{pages}{89--112}.
\newblock \URLprefix \url{https://doi.org/10.1137/S003614450036757X},
  \DOIprefix\doi{10.1137/S003614450036757X},
  \href{http://arxiv.org/abs/https://doi.org/10.1137/S003614450036757X}{{\tt
  arXiv:https://doi.org/10.1137/S003614450036757X}}.
\bibitem[{Harten and Hyman(1983)}]{harten1983}
\bibinfo{author}{Harten, A.}, \bibinfo{author}{Hyman, J.M.},
  \bibinfo{year}{1983}.
\newblock \bibinfo{title}{Self adjusting grid methods for one-dimensional
  hyperbolic conservation laws}.
\newblock \bibinfo{journal}{Journal of Computational Physics}
  \bibinfo{volume}{50}, \bibinfo{pages}{235 -- 269}.
\newblock \URLprefix
  \url{http://www.sciencedirect.com/science/article/pii/0021999183900669},
  \DOIprefix\doi{https://doi.org/10.1016/0021-9991(83)90066-9}.
\bibitem[{Hesthaven and Warburton(2008)}]{hestha1}
\bibinfo{author}{Hesthaven, J.S.}, \bibinfo{author}{Warburton, T.},
  \bibinfo{year}{2008}.
\newblock \bibinfo{title}{Nodal {Discontinuous} {Galerkin} {Methods:}
  {Algorithms,} {Analysis,} and {Applications}}.
\newblock \bibinfo{publisher}{Springer {New} {York}}.
\newblock \URLprefix \url{https://doi.org/10.1007/978-0-387-72067-8_6},
  \DOIprefix\doi{10.1007/978-0-387-72067-8_6}.
\bibitem[{Jiang and Shu(1996)}]{jiang1996}
\bibinfo{author}{Jiang, G.S.}, \bibinfo{author}{Shu, C.W.},
  \bibinfo{year}{1996}.
\newblock \bibinfo{title}{Efficient implementation of {W}eighted {ENO}
  schemes}.
\newblock \bibinfo{journal}{Journal of Computational Physics}
  \bibinfo{volume}{126}, \bibinfo{pages}{202 -- 228}.
\newblock \URLprefix
  \url{http://www.sciencedirect.com/science/article/pii/S0021999196901308},
  \DOIprefix\doi{10.1006/jcph.1996.0130}.
\bibitem[{Jin and Liu(1996)}]{jin1996}
\bibinfo{author}{Jin, S.}, \bibinfo{author}{Liu, J.G.}, \bibinfo{year}{1996}.
\newblock \bibinfo{title}{The effects of numerical viscosities: I. slowly
  moving shocks}.
\newblock \bibinfo{journal}{Journal of Computational Physics}
  \bibinfo{volume}{126}, \bibinfo{pages}{373 -- 389}.
\newblock \URLprefix
  \url{http://www.sciencedirect.com/science/article/pii/S0021999196901448},
  \DOIprefix\doi{https://doi.org/10.1006/jcph.1996.0144}.
\bibitem[{Johnsen(2013)}]{johnsen2013}
\bibinfo{author}{Johnsen, E.}, \bibinfo{year}{2013}.
\newblock \bibinfo{title}{Analysis of numerical errors generated by slowly
  moving shock waves}.
\newblock \bibinfo{journal}{AIAA Journal} \bibinfo{volume}{51},
  \bibinfo{pages}{1269--1274}.
\newblock \URLprefix \url{https://doi.org/10.2514/1.J051884},
  \DOIprefix\doi{10.2514/1.J051884},
  \href{http://arxiv.org/abs/https://doi.org/10.2514/1.J051884}{{\tt
  arXiv:https://doi.org/10.2514/1.J051884}}.
\bibitem[{Kitamura and Shima(2019)}]{kitamura2019}
\bibinfo{author}{Kitamura, K.}, \bibinfo{author}{Shima, E.},
  \bibinfo{year}{2019}.
\newblock \bibinfo{title}{Numerical experiments on anomalies from stationary,
  slowly moving, and fast-moving shocks}.
\newblock \bibinfo{journal}{AIAA Journal} \bibinfo{volume}{57},
  \bibinfo{pages}{1763--1772}.
\newblock \URLprefix \url{https://doi.org/10.2514/1.J057366},
  \DOIprefix\doi{10.2514/1.J057366},
  \href{http://arxiv.org/abs/https://doi.org/10.2514/1.J057366}{{\tt
  arXiv:https://doi.org/10.2514/1.J057366}}.
\bibitem[{Laney(1998)}]{laney1998}
\bibinfo{author}{Laney, C.B.}, \bibinfo{year}{1998}.
\newblock \bibinfo{title}{Computational gasdynamics}.
  \bibinfo{publisher}{Cambridge university press}.
  chapter~\bibinfo{chapter}{5}.
\newblock pp. \bibinfo{pages}{84--87}.
\bibitem[{Levy et~al.(1993)Levy, Powell and van Leer}]{levy1993}
\bibinfo{author}{Levy, D.W.}, \bibinfo{author}{Powell, K.G.},
  \bibinfo{author}{van Leer, B.}, \bibinfo{year}{1993}.
\newblock \bibinfo{title}{Use of a rotated riemann solver for the
  two-dimensional euler equations}.
\newblock \bibinfo{journal}{Journal of Computational Physics}
  \bibinfo{volume}{106}, \bibinfo{pages}{201 -- 214}.
\newblock \URLprefix
  \url{http://www.sciencedirect.com/science/article/pii/S0021999183711034},
  \DOIprefix\doi{https://doi.org/10.1016/S0021-9991(83)71103-4}.
\bibitem[{Lin(1995)}]{lin1995}
\bibinfo{author}{Lin, H.C.}, \bibinfo{year}{1995}.
\newblock \bibinfo{title}{Dissipation additions to flux-difference splitting}.
\newblock \bibinfo{journal}{Journal of Computational Physics}
  \bibinfo{volume}{117}, \bibinfo{pages}{20 -- 27}.
\newblock \URLprefix
  \url{http://www.sciencedirect.com/science/article/pii/S0021999185710406},
  \DOIprefix\doi{https://doi.org/10.1006/jcph.1995.1040}.
\bibitem[{Liou(2000)}]{liou2000}
\bibinfo{author}{Liou, M.S.}, \bibinfo{year}{2000}.
\newblock \bibinfo{title}{Mass flux schemes and connection to shock
  instability}.
\newblock \bibinfo{journal}{Journal of Computational Physics}
  \bibinfo{volume}{160}, \bibinfo{pages}{623 -- 648}.
\newblock \URLprefix
  \url{http://www.sciencedirect.com/science/article/pii/S0021999100964787},
  \DOIprefix\doi{https://doi.org/10.1006/jcph.2000.6478}.
\bibitem[{Liou(2006)}]{liou2006}
\bibinfo{author}{Liou, M.S.}, \bibinfo{year}{2006}.
\newblock \bibinfo{title}{A sequel to ausm, part ii: Ausm+-up for all speeds}.
\newblock \bibinfo{journal}{Journal of Computational Physics}
  \bibinfo{volume}{214}, \bibinfo{pages}{137 -- 170}.
\newblock \URLprefix
  \url{http://www.sciencedirect.com/science/article/pii/S0021999105004274},
  \DOIprefix\doi{https://doi.org/10.1016/j.jcp.2005.09.020}.
\bibitem[{Liu et~al.(1994)Liu, Osher and Chan}]{liu1994}
\bibinfo{author}{Liu, X.D.}, \bibinfo{author}{Osher, S.},
  \bibinfo{author}{Chan, T.}, \bibinfo{year}{1994}.
\newblock \bibinfo{title}{Weighted essentially non-oscillatory schemes}.
\newblock \bibinfo{journal}{Journal of Computational Physics}
  \bibinfo{volume}{115}, \bibinfo{pages}{200 -- 212}.
\newblock \URLprefix
  \url{http://www.sciencedirect.com/science/article/pii/S0021999184711879},
  \DOIprefix\doi{10.1006/jcph.1994.1187}.
\bibitem[{Merriman(2003)}]{merriman2003}
\bibinfo{author}{Merriman, B.}, \bibinfo{year}{2003}.
\newblock \bibinfo{title}{Understanding the {S}hu-{O}sher conservative finite
  difference form}.
\newblock \bibinfo{journal}{Journal of Scientific Computing}
  \bibinfo{volume}{19}, \bibinfo{pages}{309--322}.
\newblock \URLprefix \url{https://doi.org/10.1023/A:1025312210724},
  \DOIprefix\doi{10.1023/A:1025312210724}.
\bibitem[{Nishikawa and Kitamura(2008)}]{nishikawa2008}
\bibinfo{author}{Nishikawa, H.}, \bibinfo{author}{Kitamura, K.},
  \bibinfo{year}{2008}.
\newblock \bibinfo{title}{Very simple, carbuncle-free,
  boundary-layer-resolving, rotated-hybrid riemann solvers}.
\newblock \bibinfo{journal}{Journal of Computational Physics}
  \bibinfo{volume}{227}, \bibinfo{pages}{2560 -- 2581}.
\newblock \URLprefix
  \url{http://www.sciencedirect.com/science/article/pii/S0021999107004822},
  \DOIprefix\doi{https://doi.org/10.1016/j.jcp.2007.11.003}.
\bibitem[{Peng et~al.(2019)Peng, Zhai, Ni, Yong and Shen}]{peng2019}
\bibinfo{author}{Peng, J.}, \bibinfo{author}{Zhai, C.}, \bibinfo{author}{Ni,
  G.}, \bibinfo{author}{Yong, H.}, \bibinfo{author}{Shen, Y.},
  \bibinfo{year}{2019}.
\newblock \bibinfo{title}{An adaptive characteristic-wise reconstruction weno-z
  scheme for gas dynamic euler equations}.
\newblock \bibinfo{journal}{Computers \& Fluids} \bibinfo{volume}{179},
  \bibinfo{pages}{34 -- 51}.
\newblock \URLprefix
  \url{http://www.sciencedirect.com/science/article/pii/S0045793018304997},
  \DOIprefix\doi{https://doi.org/10.1016/j.compfluid.2018.08.008}.
\bibitem[{Qiu et~al.(2006)Qiu, Koo and Shu}]{qks}
\bibinfo{author}{Qiu, J.}, \bibinfo{author}{Koo, B.C.}, \bibinfo{author}{Shu,
  C.W.}, \bibinfo{year}{2006}.
\newblock \bibinfo{title}{A numerical study for the performance of the
  {Runge-Kutta} discontinuous {Galerkin} method based on different numerical
  fluxes.}
\newblock \bibinfo{journal}{Journal of Computational Physics}
  \bibinfo{volume}{212}, \bibinfo{pages}{540--565}.
\bibitem[{Qiu and Shu(2002)}]{qiu2002}
\bibinfo{author}{Qiu, J.}, \bibinfo{author}{Shu, C.W.}, \bibinfo{year}{2002}.
\newblock \bibinfo{title}{On the construction, comparison, and local
  characteristic decomposition for high-order central weno schemes}.
\newblock \bibinfo{journal}{Journal of Computational Physics}
  \bibinfo{volume}{183}, \bibinfo{pages}{187 -- 209}.
\newblock \URLprefix
  \url{http://www.sciencedirect.com/science/article/pii/S0021999102971913},
  \DOIprefix\doi{https://doi.org/10.1006/jcph.2002.7191}.
\bibitem[{Qiu and Shu(2005)}]{qs2}
\bibinfo{author}{Qiu, J.}, \bibinfo{author}{Shu, C.W.}, \bibinfo{year}{2005}.
\newblock \bibinfo{title}{A comparison of troubled-cell indicators for
  {Runge-Kutta} discontinuous {Galerkin} methods using weighted essentially
  nonoscillatory limiters.}
\newblock \bibinfo{journal}{SIAM J. Sci. Comput.} \bibinfo{volume}{27},
  \bibinfo{pages}{995--1013}.
\bibitem[{Reed and Hill(1973)}]{rh}
\bibinfo{author}{Reed, W.}, \bibinfo{author}{Hill, T.}, \bibinfo{year}{1973}.
\newblock \bibinfo{title}{Triangular mesh methods for the neutron transport
  equation.}
\newblock \bibinfo{type}{Technical Report} \bibinfo{number}{Technical Report
  LA-UR-73-479}. Los Alamos Scientific Laboratory.
\bibitem[{Roberts(1990)}]{roberts1990}
\bibinfo{author}{Roberts, T.W.}, \bibinfo{year}{1990}.
\newblock \bibinfo{title}{The behavior of flux difference splitting schemes
  near slowly moving shock waves}.
\newblock \bibinfo{journal}{Journal of Computational Physics}
  \bibinfo{volume}{90}, \bibinfo{pages}{141 -- 160}.
\newblock \URLprefix
  \url{http://www.sciencedirect.com/science/article/pii/002199919090200K},
  \DOIprefix\doi{https://doi.org/10.1016/0021-9991(90)90200-K}.
\bibitem[{Roe(1981)}]{roe1981}
\bibinfo{author}{Roe, P.}, \bibinfo{year}{1981}.
\newblock \bibinfo{title}{Approximate riemann solvers, parameter vectors, and
  difference schemes}.
\newblock \bibinfo{journal}{Journal of Computational Physics}
  \bibinfo{volume}{43}, \bibinfo{pages}{357 -- 372}.
\newblock \URLprefix
  \url{http://www.sciencedirect.com/science/article/pii/0021999181901285},
  \DOIprefix\doi{https://doi.org/10.1016/0021-9991(81)90128-5}.
\bibitem[{Shu(1988)}]{shu}
\bibinfo{author}{Shu, C.W.}, \bibinfo{year}{1988}.
\newblock \bibinfo{title}{{TVD} time discretizations.}
\newblock \bibinfo{journal}{SIAM J. Sci. Stat. Comput.} \bibinfo{volume}{9},
  \bibinfo{pages}{1073--1084}.
\bibitem[{Shu(2009)}]{shu2009}
\bibinfo{author}{Shu, C.W.}, \bibinfo{year}{2009}.
\newblock \bibinfo{title}{High order weighted essentially nonoscillatory
  schemes for convection dominated problems}.
\newblock \bibinfo{journal}{SIAM Review} \bibinfo{volume}{51},
  \bibinfo{pages}{82--126}.
\newblock \URLprefix \url{https://doi.org/10.1137/070679065},
  \DOIprefix\doi{10.1137/070679065},
  \href{http://arxiv.org/abs/https://doi.org/10.1137/070679065}{{\tt
  arXiv:https://doi.org/10.1137/070679065}}.
\bibitem[{Spiegel et~al.(2015)Spiegel, Huynh and DeBonis}]{spiegel2015}
\bibinfo{author}{Spiegel, S.C.}, \bibinfo{author}{Huynh, H.},
  \bibinfo{author}{DeBonis, J.R.}, \bibinfo{year}{2015}.
\newblock \bibinfo{title}{A survey of the isentropic euler vortex problem using
  high-order methods}, in: \bibinfo{booktitle}{22nd AIAA Computational Fluid
  Dynamics Conference}, p. \bibinfo{pages}{2444}.
\bibitem[{Steger and Warming(1981)}]{steger1981}
\bibinfo{author}{Steger, J.L.}, \bibinfo{author}{Warming, R.},
  \bibinfo{year}{1981}.
\newblock \bibinfo{title}{Flux vector splitting of the inviscid gasdynamic
  equations with application to finite-difference methods}.
\newblock \bibinfo{journal}{Journal of Computational Physics}
  \bibinfo{volume}{40}, \bibinfo{pages}{263 -- 293}.
\newblock \URLprefix
  \url{http://www.sciencedirect.com/science/article/pii/0021999181902102},
  \DOIprefix\doi{https://doi.org/10.1016/0021-9991(81)90210-2}.
\bibitem[{Stiriba and Donat(2003)}]{stiriba2003}
\bibinfo{author}{Stiriba, Y.}, \bibinfo{author}{Donat, R.},
  \bibinfo{year}{2003}.
\newblock \bibinfo{title}{A numerical study of postshock oscillations in slowly
  moving shock waves}.
\newblock \bibinfo{journal}{Computers \& Mathematics with Applications}
  \bibinfo{volume}{46}, \bibinfo{pages}{719 -- 739}.
\newblock \URLprefix
  \url{http://www.sciencedirect.com/science/article/pii/S0898122103901374},
  \DOIprefix\doi{https://doi.org/10.1016/S0898-1221(03)90137-4}.
\bibitem[{Toro(2013)}]{toro2013}
\bibinfo{author}{Toro, E.F.}, \bibinfo{year}{2013}.
\newblock \bibinfo{title}{Riemann solvers and numerical methods for fluid
  dynamics: a practical introduction}.
\newblock \bibinfo{publisher}{Springer Science \& Business Media}.
\bibitem[{Wu et~al.(2016)Wu, Liang and Zhao}]{wu2016}
\bibinfo{author}{Wu, X.}, \bibinfo{author}{Liang, J.}, \bibinfo{author}{Zhao,
  Y.}, \bibinfo{year}{2016}.
\newblock \bibinfo{title}{A new smoothness indicator for third-order weno
  scheme}.
\newblock \bibinfo{journal}{International Journal for Numerical Methods in
  Fluids} \bibinfo{volume}{81}, \bibinfo{pages}{451--459}.
\newblock \URLprefix
  \url{https://onlinelibrary.wiley.com/doi/abs/10.1002/fld.4194},
  \DOIprefix\doi{10.1002/fld.4194},
  \href{http://arxiv.org/abs/https://onlinelibrary.wiley.com/doi/pdf/10.1002/fld.4194}{{\tt
  arXiv:https://onlinelibrary.wiley.com/doi/pdf/10.1002/fld.4194}}.
\bibitem[{Xu(1999)}]{xu1999}
\bibinfo{author}{Xu, K.}, \bibinfo{year}{1999}.
\newblock \bibinfo{title}{Does perfect Riemann solver exist?}.
  \bibinfo{publisher}{14th Computational Fluid Dynamics Conference, AIAA}.
\newblock pp. \bibinfo{pages}{771--779}.
\newblock \URLprefix \url{https://arc.aiaa.org/doi/abs/10.2514/6.1999-3344},
  \DOIprefix\doi{10.2514/6.1999-3344},
  \href{http://arxiv.org/abs/https://arc.aiaa.org/doi/pdf/10.2514/6.1999-3344}{{\tt
  arXiv:https://arc.aiaa.org/doi/pdf/10.2514/6.1999-3344}}.
\bibitem[{Zhang et~al.(2011)Zhang, Jiang and Shu}]{zhang2011}
\bibinfo{author}{Zhang, S.}, \bibinfo{author}{Jiang, S.}, \bibinfo{author}{Shu,
  C.W.}, \bibinfo{year}{2011}.
\newblock \bibinfo{title}{Improvement of convergence to steady state solutions
  of~ {E}uler equations with~ the~ {WENO}~ schemes}.
\newblock \bibinfo{journal}{Journal of Scientific Computing}
  \bibinfo{volume}{47}, \bibinfo{pages}{216--238}.
\newblock \URLprefix \url{https://doi.org/10.1007/s10915-010-9435-5},
  \DOIprefix\doi{10.1007/s10915-010-9435-5}.
\bibitem[{Zhong and Shu(2013)}]{zs}
\bibinfo{author}{Zhong, X.}, \bibinfo{author}{Shu, C.W.}, \bibinfo{year}{2013}.
\newblock \bibinfo{title}{A simple weighted essentially nonoscillatory limiter
  for {Runge-Kutta} discontinuous {Galerkin} methods.}
\newblock \bibinfo{journal}{Journal of Computational Physics}
  \bibinfo{volume}{232}, \bibinfo{pages}{397--415}.

\end{thebibliography}
\bibliographystyle{elsarticle/elsarticle-harv}

\end{document}